\documentclass[12pt]{article}
\usepackage[left=1.1in,right=1.1in,top=1in,bottom=1in]{geometry}
\usepackage[T1]{fontenc}
\usepackage{lmodern,microtype}
\usepackage{amsmath,amssymb,amsthm,mathtools}
\usepackage{enumitem}
\usepackage{tikz}
\usetikzlibrary{calc,decorations.pathreplacing,arrows.meta,shapes.geometric}
\usepackage{caption}
\usepackage{subcaption}
\usepackage[section]{placeins}
\usepackage{needspace}
\definecolor{prodA}{HTML}{245C86}
\definecolor{prodB}{HTML}{A33C50}
\definecolor{prodC}{HTML}{287765}
\definecolor{prodZ}{HTML}{9B6418}
\tikzset{
  product/.style={x=0.67cm,y=0.53cm,font=\footnotesize,line cap=round,line join=round},
  prod guide/.style={draw=black!17,line width=0.3pt,densely dotted},
  prod A/.style={circle,draw=prodA,fill=prodA,inner sep=0pt,minimum size=4.6pt},
  prod B/.style={rectangle,draw=prodB,fill=prodB,inner sep=0pt,minimum size=4.6pt},
  prod C/.style={diamond,draw=prodC,fill=prodC,inner sep=0pt,minimum size=5.8pt},
  prod Z/.style={circle,draw=prodZ,line width=1pt,fill=white,inner sep=0pt,minimum size=5.5pt},
  prod edge/.style={draw=black!75,line width=0.7pt},
  prod focus/.style={draw=prodC,line width=1.0pt},
  prod muted/.style={draw=black!38,line width=0.6pt,dashed},
  prod brace/.style={decorate,decoration={brace,amplitude=3pt},draw=black!65,line width=0.4pt},
  prod note/.style={font=\scriptsize,align=center,inner sep=2pt}
}
\newcommand{\ProdGrid}[2]{%
  \foreach \xx in {1,...,#1}{\draw[prod guide] (\xx,0.7)--(\xx,{#2+0.3});}%
  \foreach \yy in {1,...,#2}{\draw[prod guide] (0.7,\yy)--({#1+0.3},\yy);}%
}
\newcommand{\ProdX}[2]{\node[anchor=north,inner sep=2pt] at (#1,0.58) {$#2$};}
\newcommand{\ProdY}[2]{\node[anchor=east,inner sep=2pt] at (0.55,#1) {$#2$};}
\newcommand{\ProdPoint}[3]{\node[prod #3] at (#1,#2) {};}
\newcommand{\ProdPair}[4]{%
  \draw[prod brace] (#1,{#3+0.5})--(#2,{#3+0.5})
    node[midway,above=4pt] {$#4$};}

\usepackage[hidelinks]{hyperref}
\numberwithin{equation}{section}
\newtheorem{theorem}{Theorem}[section]
\newtheorem{proposition}[theorem]{Proposition}
\newtheorem{lemma}[theorem]{Lemma}
\newtheorem{corollary}[theorem]{Corollary}
\newtheorem{claim}[theorem]{Claim}
\theoremstyle{definition}
\newtheorem{definition}[theorem]{Definition}
\theoremstyle{remark}
\newtheorem{observation}[theorem]{Observation}

\newcommand{\F}{\mathbb F_2}
\newcommand{\U}{\mathcal U}
\newcommand{\col}{\operatorname{col}}
\newcommand{\sgn}{\operatorname{sgn}}
\newcommand{\spn}{\operatorname{span}}
\newcommand{\odd}{\textnormal{-odd}}
\newcommand{\uniq}{\textnormal{-unique}}
\newcommand{\ff}{\mathbf f}
\newcommand{\pp}{\mathbf p}
\newcommand{\ttv}{\mathbf t}
\newcommand{\ee}{\mathbf e}
\newcommand{\hh}{\mathbf h}

\title{A variant of the Erd\H{o}s--Gy\'arf\'as problem for $K_9$ and $K_{12}$
}
\author{Luyi Li \footnote{Research Center for Applied Mathematics, Renmin University of China, Beijing, China. 
{\tt Email: liluyiplus@gmail.com} }, \ \
Guiying Yan \footnote{Academy of Mathematics and Systems Science, Chinese Academy of Sciences, Beijing,  China. 
{\tt Email: yangy@amt.ac.cn} }, \ \
and Zicheng Ye\footnote{Corresponding author: Academy of Mathematics and Systems Science, Chinese Academy of Sciences, Beijing,  China.  {\tt
Email: yezicheng@amss.ac.cn}}
\thanks{This work was supported by the National Key R\&D Program of China (2023YFA1009602).}
}

\date{}
 
\begin{document}
\maketitle
\begin{abstract}
For a graph $H$, let $r_H(n)$ be the minimum number of colors in an edge-coloring of $K_n$ such that every copy of $H$ contains a color appearing an odd number of times. We prove that $r_{K_9}(n)\le\exp(O((\log n)^{2/3}))$ and $r_{K_{12}}(n)\le\exp(O((\log n)^{2/3}))$. These bounds on $r_H(n)$ yield linear graph codes of density at least $\exp(-O((\log n)^{2/3}))$ that avoid images of $K_9$ and $K_{12}$.
\end{abstract}

\section{Introduction}

The Erd\H{o}s--Gy\'arf\'as problem asks for the minimum number of colors in an edge-coloring of $K_n$ in which every copy of $K_p$ receives at least $q$ colors; see~\cite{Erdos1997variant}. Motivated by graph codes, Alon~\cite{Alon2024European} considered a related condition on the multiplicities of the colors. For a graph $H$, an edge-coloring is called \emph{$H$-odd} if every copy of $H$ contains a color appearing an odd number of times. We write $r_H(n)$ for the least number of colors in an $H$-odd edge-coloring of $K_n$. A copy of $H$ is \emph{even-chromatic} if each of its colors appears an even number of times. Thus an $H$-odd coloring is precisely a coloring with no even-chromatic copy of $H$.

If $H$ has an odd number of edges, every edge-coloring is $H$-odd. For complete graphs, the nontrivial cases are therefore $K_t$ with $t\equiv0,1\pmod4$. Ge, Xu and Zhang~\cite{Ge2023new} conjectured that
\begin{equation}\label{eq:conjecture}
 r_{K_t}(n)=n^{o(1)} \qquad\text{for every fixed }t\ge4\text{ with }t\equiv0,1\pmod4.
\end{equation}
The case $K_4$ follows from the construction of Cameron and Heath~\cite{Cameron2023New}. Ge, Xu and Zhang~\cite{Ge2023new} proved the case $K_5$ by adapting constructions for generalized Ramsey problems; Bennett, Heath and Zerbib~\cite{Bennett2023Edge} obtained this result independently. Yip~\cite{Yip2026variant} introduced an amalgamation of two edge-colorings and proved that
$$
r_{K_8}(n)\le \exp\bigl(O((\log n)^{2/3})\bigr).
$$
The same amalgamation method gives colorings using $\exp(O(\sqrt{\log n}))$ colors in which every copy of $K_4$ and every copy of $K_5$ contains a color appearing exactly once.

We extend the amalgamation argument to $K_9$ and $K_{12}$. During the construction, we retain a stronger uniqueness condition. We call a coloring \emph{$H$-unique} if every copy of $H$ contains a color appearing on exactly one edge. Clearly $H$-uniqueness implies $H$-oddness. We obtain the following two bounds.

\begin{theorem}\label{thm:k9}
There is an absolute constant $C_9>0$ such that, for every $n\ge 9$,
$$
r_{K_9}(n)\le \exp\bigl(C_9(\log n)^{2/3}\bigr).
$$
In particular, $r_{K_9}(n)=n^{o(1)}$.
\end{theorem}

\begin{theorem}\label{thm:k12}
There is an absolute constant $C_{12}>0$ such that, for every $n\ge 12$,
$$
r_{K_{12}}(n)\le \exp\bigl(C_{12}(\log n)^{2/3}\bigr).
$$
In particular, $r_{K_{12}}(n)=n^{o(1)}$.
\end{theorem}

The proofs use the same construction and color-count estimate as in~\cite{Yip2026variant}, but proceed by separate inductions. For Theorem~\ref{thm:k9}, we preserve $K_3$-uniqueness, $K_4$-uniqueness and $K_5$-uniqueness together with $K_9$-oddness. For Theorem~\ref{thm:k12}, we preserve $K_3$-uniqueness, $K_4$-uniqueness and $K_6$-uniqueness together with $K_8$-oddness and $K_{12}$-oddness. In each local preservation argument, the exact vertical colors and a projection onto the $\mathcal{C}$-color first restrict the column supports of an assumed even-chromatic set. For nine vertices, this leaves a doubled four-set and one vertex, or two rectangles and one vertex. For twelve vertices, a doubled six-set is excluded by an immediate diagonal count, leaving a doubled four-set and a rectangle, or three rectangles. We then exclude the remaining configurations using parity equations for blocks of edges. The additional  vertex in the nine-vertex case requires us to account for coincident rows. In the twelve- vertex case, after simpler incidences have been excluded, the three row pairs form either a path of length three or a path of length two together with a disjoint edge. We give explicit proofs for both patterns.

We also derive a consequence for graph codes. Identify graphs on $[n]$ with their edge-indicator vectors in $\F^m$ with $m={\binom n2}$, so that vector addition corresponds to symmetric difference. An \emph{image} of $H$ on $[n]$ is a graph consisting of a copy of $H$ and isolated vertices. A \emph{$H$-code} is a set $\mathcal{C}$ of graphs such that no difference of two distinct members is an image of a member of $H$. This is the graph-code convention of Alon et al.~\cite{Alon2023Structured}. If $\mathcal{C}$ is a vector subspace, then $\mathcal{C}$ is a \emph{linear graph code}. The relation between $H$-odd edge-colorings and linear graph codes $\mathcal{C}$ was developed in~\cite{Alon2024European,Versteegen2025Random}. It was shown that $r_H(n)^{-1}\leq |\mathcal{C}|/{2^{\binom n2}} $. Using Theorem \ref{thm:k9} and \ref{thm:k12} and this relation, we obtain the following corollary.

\begin{corollary}\label{cor:code}
Let $H=K_9$ or $K_{12}$. There is an absolute constant $C_H>0$ such that, for every $n\ge2$, there is a linear $H$-code $\mathcal C \subseteq \F^{\binom n2}$ satisfying
$$
\frac{|\mathcal C|}{2^{\binom n2}}
 \ge \exp\bigl(-C_H(\log n)^{2/3}\bigr).
$$
\end{corollary}

Section~\ref{sec:construction} introduces the amalgamation, constructs the auxiliary colorings and gives a common color-count estimate. Section~\ref{sec:parity} develops the common parity notation and the two-rectangle argument. Sections~\ref{sec:nine} and \ref{sec:twelve} give the local preservation arguments and complete the proofs of Theorems~\ref{thm:k9} and~\ref{thm:k12}, respectively. The latter section also proves the preservation of $K_6$-uniqueness. 

Throughout, all graphs are finite and simple. We write $[n]=\{1,\ldots,n\}$, and $\binom X2$ for the set of unordered pairs from $X$. For a coloring $c$, the notation $c(xx')$ refers to the color of $\{x,x'\}$; if $A$ is a two-element set, we also write $c(A)$. All parity vectors are over $\F$. All logarithms are natural unless a base is specified. A condition on $K_s$ is understood vacuously on sets with fewer than $s$ vertices. Labels in small configurations specify incidence and impose no order on the rows or columns.

\section{The amalgamation coloring construction}\label{sec:construction}

We use Yip's amalgamation~\cite{Yip2026variant}, expressed through three disjoint palettes rather than four color coordinates. The two descriptions induce the same equivalence relation on edges, since the exact row pair of a vertical edge already determines its $\mathcal{D}$-color.

\begin{definition}[\cite{Yip2026variant}]\label{def:amalgamation}
Let $c$ and $d$ be edge-colorings of complete graphs on ordered sets $X$ and $Y$, respectively. Their \emph{amalgamation} $c\otimes d$ is an edge-coloring of the complete graph on $X\times Y$. For distinct
$(x,y)$ and $(x',y')$, put
\begin{equation*}
 (c\otimes d)((x,y),(x',y'))=
 \begin{cases}
 (V,\{y,y'\}),&x=x',\\
 (H,c(xx')),&y=y',\\
 (D,c(xx'),d(yy'),\sgn((x'-x)(y'-y))),&x\ne x',\ y\ne y'.
 \end{cases}
\end{equation*}
The colorings on $X$ and $Y$ are denoted as denoted as $\mathcal{C}$-color and $\mathcal{D}$-color, respectively. The symbols $V,H,D$ distinguish vertical, horizontal and diagonal edges. The sign in the last line is positive when the two coordinate orders agree, and negative otherwise. See Fig.~\ref{fig:amalgamation} for an example.
\end{definition}

In all figures, $X$ is placed horizontally and $Y$ vertically, so vertices with the same first coordinate lie in a column. Dotted lines are coordinate guides rather than graph edges. The vertex shapes and display colors distinguish product parts and do not specify values of $c$ or $d$. A hollow vertex denotes the additional singleton. Each panel displays its full vertex set, but only explicitly indicated graph edges are drawn. Unless an order is specified in a caption, coordinate positions show incidence only and impose no order assumption.

\begin{figure}[!htbp]
\centering
\begin{minipage}[c]{0.43\linewidth}
\centering
\begin{tikzpicture}[product,x=1.75cm,y=1.35cm]
\ProdGrid{2}{2}
\ProdX{1}{x}\ProdX{2}{x'}
\ProdY{1}{y}\ProdY{2}{y'}
\draw[draw=prodA,line width=1pt] (1,1)--(1,2) (2,1)--(2,2);
\draw[draw=prodB,line width=1pt] (1,1)--(2,1) (1,2)--(2,2);
\draw[draw=prodC,line width=1pt] (1,1)--(2,2);
\draw[draw=prodZ,line width=1pt,dashed] (1,2)--(2,1);
\node[fill=white,inner sep=1pt,text=prodC] at (1.29,1.29) {$+$};
\node[fill=white,inner sep=1pt,text=prodZ] at (1.72,1.28) {$-$};
\foreach \xx in {1,2}{\foreach \yy in {1,2}{\fill (\xx,\yy) circle[radius=2.2pt];}}
\end{tikzpicture}
\end{minipage}\hfill
\begin{minipage}[c]{0.54\linewidth}
\small
\renewcommand{\arraystretch}{1.45}
\begin{tabular}{@{}ll@{}}
\textcolor{prodA}{Vertical} & $(V,\{y,y'\})$\\
\textcolor{prodB}{Horizontal} & $(H,c(xx'))$\\
\textcolor{prodC}{Positive diagonal} & $(D,c(xx'),d(yy'),+)$\\
\textcolor{prodZ}{Negative diagonal} & $(D,c(xx'),d(yy'),-)$
\end{tabular}
\end{minipage}
\caption{The amalgamation in Definition~\ref{def:amalgamation}, illustrated on $\{x,x'\}\times\{y,y'\}$ with $x<x'$ and $y<y'$. The two vertical edges have the same color, and the two horizontal edges have the same color. The diagonal signs distinguish the two remaining edges.}
\label{fig:amalgamation}
\end{figure} 

For $x\ne x'$ and $y\ne y'$, the two diagonal edges $(x,y)(x',y')$ and $(x,y')(x',y)$ receive different colors, since their signs are opposite. We call $\{x\}\times Y$ a \emph{column} and $X\times\{y\}$ a \emph{row}. For a coloring $c$, let $\col(c)$ denote its number of colors.

\begin{lemma}\label{lem:count}
If $\col(c) = a$, $\col(d) = b$ and $|Y|=m$, then
$\col(c\otimes d)\le (2b+1)a+\binom m2.$
\end{lemma}
\begin{proof}
The vertical palette has at most $\binom m2$ colors. Horizontal edges use at most $a$ colors. A diagonal color is determined by a color of $c$, a color of $d$ and one of two signs, so there are at most $2ab$ diagonal colors. These palettes are disjoint.
\end{proof}

We next recall preservation properties of the amalgamation from~\cite{Yip2026variant}. These will be used to construct the auxiliary colorings needed later. 

\begin{proposition}[\cite{Yip2026variant}]
\label{prop:unique}
For any edge-coloring $d$, the following hold.
\begin{enumerate}[label=\textup{(\roman*)},leftmargin=2em]
\item For $s\in\{3,4\}$, if $c$ is $K_s$-unique, then $c\otimes d$ is $K_s$-unique.
\item If $c$ is both $K_3$-unique and $K_5$-unique, then $c\otimes d$ has both properties.
\end{enumerate}
\end{proposition}

Let $\U$ denote the simultaneous property
\begin{equation*}
 \U=\{K_3\uniq,\ K_4\uniq,\ K_5\uniq\}.
\end{equation*}
By Proposition~\ref{prop:unique}, amalgamation with an arbitrary edge-coloring $d$ preserves $\U$.

\begin{proposition}\label{prop:aux}
For every $n\ge2$, there is a coloring with property $\U$ using at most $\exp(C_0\sqrt{\log n})$ colors, where $C_0$ is an absolute constant.
\end{proposition}
\begin{proof}
Fix an integer $m\ge2$, and let $d$ be a monochromatic coloring on $m$ vertices. Starting with the edge-free coloring on one vertex, repeatedly amalgamate with $d$. After $j$ steps, the resulting coloring has $m^j$ vertices and retains property $\U$ by Proposition~\ref{prop:unique}. If $a_j$ denotes its number of colors, Lemma~\ref{lem:count} gives
$$
 a_0=0 \text{ and }  a_{j+1}\le3a_j+\binom m2.
$$
Consequently,
\begin{equation*}
 a_j\le\binom m2\frac{3^j-1}{2}\le m^2 3^j.
\end{equation*}
For $n\ge3$, take $k=\lceil\sqrt{\log n}\rceil$ and $m=\lceil n^{1/k}\rceil$. Then $m^k\ge n$ and
$$
 \log m\le\frac{\log n}{k}+\log2=O(\sqrt{\log n}),
$$
while $k=O(\sqrt{\log n})$. The bound $a_k\le m^2 3^k$ therefore gives the claimed estimate on $m^k$ vertices. Restricting this coloring to any $n$ vertices preserves property  and does not increase the number of colors. The case $n=2$ follows by increasing $C_0$ if necessary.
\end{proof}

The following lemma is an application of Lemma 2.5 in~\cite{Yip2026variant} with $p=2/3$ and $q=1/2$.

\begin{lemma}[\cite{Yip2026variant}]\label{lem:iteration}
Let $\mathcal P$ be a property of edge-colorings on finite ordered sets that holds on a single vertex and is preserved under restriction. Suppose that $c\otimes d$ has property $\mathcal P$ whenever $c$ has property $\mathcal P$ and $d$ has property $\U$. Then, for every $n\ge2$, there is a coloring of $K_n$ with property $\mathcal P$ using at most $\exp(C(\log n)^{2/3})$ colors, where $C$ is an absolute constant.
\end{lemma}

\section{Auxiliary lemmas}\label{sec:parity}
For a set $S\subseteq X\times Y$, set
$
 R_x=\{y:(x,y)\in S\}, m_x=|R_x|\mbox{ and }
 X_S=\{x:m_x>0\}.
$
Merging colors preserves even-chromaticity. We repeatedly use such projections to extract consequences of the parity conditions.

\begin{lemma}\label{lem:projection}
If $S$ is even-chromatic under $c\otimes d$, then the following hold.
\begin{enumerate}[label=\textup{(\roman*)},leftmargin=2em]
\item The columns of odd size induce an even-chromatic clique under $c$.
\item For  two distinct rows $u,v\in Y$, 
$|\{x\in X: \{u,v\}\subseteq R_x\}|$ is even.
\item If $|R_x|\geq 2$ for some $x\in X$, then $|R_x|\leq |S|/2$.
\end{enumerate}
\end{lemma}
\begin{proof}
(i): Fix a $\mathcal{C}$-color $\alpha$. Project every nonvertical edge to its
$\mathcal{C}$-coordinate by
$\pi_{\mathcal C}(H,\alpha)=\alpha \mbox{ and  }
\pi_{\mathcal C}(D,\alpha,\gamma,\pm)=\alpha.$
Set
$N_\alpha=|\{e\in E(S): e\text{ is nonvertical and }\pi_{\mathcal C}(e)=\alpha\}|.$
Then
$$col(N_{\alpha})=\{(H,\alpha), (D,\alpha,\gamma,+), (D,\alpha,\gamma,-)\}.$$
Let $m_S(\kappa)$ be the multiplicity of $\kappa\in \col(N_{\alpha})$ in $S$. 
Since $S$ is even-chromatic,  $m_S(\kappa)$ is even.
Hence, 
$$N_\alpha=\sum_{\kappa\in col(N_{\alpha})} m_S(\kappa)\equiv 0
\pmod 2.$$
On the other hand, every nonvertical edge joins a unique pair of distinct
columns $\{x,x'\}$. If $c(xx')=\alpha$, then $\pi_{\mathcal{C}}(e)=\alpha$ for each edge $e$
 between two columns $x$ and $x'$. Therefore
\[
N_\alpha=\sum_{\{x,x'\}\subseteq X,\; c(xx')=\alpha} m_xm_{x'}.
\] 
Note that the product $m_xm_{x'}$ is odd precisely when both
$m_x$ and $m_{x'}$ are odd. Let
$O=\{x\in X: m_x\text{ is odd}\}$
be the set of odd-sized columns. Then
\[
0\equiv N_\alpha\equiv  \big|\big\{\{x,x'\}\subseteq O: c(xx')=\alpha\big\}\big|
\pmod 2.
\]
Thus color $\alpha$ appears an even number of times in the complete graph
$K_O$ on vertex set $O$ with edge coloring $c(xx')$. Since $\alpha$ was
arbitrary, $K_O$ is an even-chromatic clique under $c$.

(ii):  Note that the color of vertical edge is 
$(V,\{u,v\})$ for any two distinct rows $u,v$. Hence
\[
\big\{e\in E(S):\operatorname{col}(e)=(V,\{u,v\})\big\}
=\big\{\{(x,u),(x,v)\}:x\in X,\;\{u,v\}\subseteq R_x\big\}.
\]
Therefore the multiplicity of $(V,\{u,v\})$ is exactly
\[
\big|\{x\in X:\{u,v\}\subseteq R_x\}\big|.
\]
Since $S$ is even-chromatic, this multiplicity must be even.

(iii): 
Fix a column $x$ with $|R_x|\geq 2$. 
Choose a row $u\in R_x$, there is another row $v\in R_x\setminus \{u\}$.
By (ii), the number of columns
containing $\{u,v\}$ is even.
Thus,
there exists another column $x'\neq x$ such that $\{u,v\}\subseteq R_{x'}$. 
It means $(x',u)\in S$.
Hence, distinct choices of $u$ in $R_x$ give distinct vertices outside column $x$.
Note that there are $|R_x|$ vertices in the column $x$, and so $|R_x|\le|S|/2$.
\end{proof}

Let $W$ be the $\mathbb{F}_2$-vector space with basis
$\{\mathbf{h}\}\cup\{\mathbf{e}_{g,+},\mathbf{e}_{g,-}:g\text{ is a color of }d\}.$
Let $\tau$ be the involution: $\tau(\hh)=\hh$, $\tau(\ee_{g,+})=\ee_{g,-}$ and $\tau(\ee_{g,-})=\ee_{g,+}$.
For row sets $U,V$, define the \emph{block signature}
\begin{equation}\label{eq:signature}
 \ff(U,V)=\sum_{u\in U,\ v\in V}
 \begin{cases}
 \hh,&u=v,\\
 \ee_{d(uv),\sgn(v-u)},&u\ne v.
 \end{cases}
\end{equation} 
Then $\ff(V,U)=\tau(\ff(U,V))$. In particular, reversing the order of the column supports preserves a zero block signature. For each $\mathcal{C}$-color $\alpha$, evenness of the nonvertical colors is equivalent to
\begin{equation}\label{eq:parityequation}
 \sum_{\substack{x<x'\text{ in }X_S, \mbox{ and }c(xx')=\alpha}}
 \ff(R_x,R_{x'})=0.
\end{equation}
Then we could obtain the following result:
\begin{observation}\label{obs:uniquezero}
If a $\mathcal{C}$-color occurs on exactly one edge $xx'$ of $K_{X_S}$, then its block signature $\ff(R_x,R_{x'})$ is zero whenever $S$ is even-chromatic.
\end{observation}

In the arguments below, we discard a set of edges only after showing that every full color occurs an even number of times within it; its contribution can then be removed from the parity equations. The remaining edges need not form an induced complete graph. If all edges incident to a vertex are discarded, however, the remaining equations are those for the set obtained by deleting that vertex.

Let $E$ be the vector space with basis $\{\mathbf{e}_g : g \text{ is a color of } d\}.$
 Define $\pi:W\to E$ by
$$
 \pi(\hh)=0\mbox{ and }\pi(\ee_{g,+})=\pi(\ee_{g,-})=\ee_g.
$$
Set $\ttv_g=\ee_{g,+}+\ee_{g,-}$ and $\ttv_v=\sum_g v_g\ttv_g$ for $v=\sum_g v_g\ee_g\in E$. The linear map $v\mapsto\ttv_v$ from $E$ to $W$ is injective. For every $w\in W$,
\begin{equation}\label{eq:taupi}
 w+\tau w=\ttv_{\pi w}.
\end{equation}
For a two-element row set $P=\{u,v\}$ we have
\begin{equation}\label{eq:rectanglevector}
 \ff(P,P)=\ttv_{d(P)}\ne0.
\end{equation}
We abbreviate $\ttv_{d(P)}$ to $\ttv_P$. The coefficient of $\hh$ in $\ff(U,V)$ is $|U\cap V|$ modulo two. The map $\pi\ff$ is bilinear in row sets under symmetric difference. Repeated ordered row pairs in~\eqref{eq:signature} cancel, and terms with equal rows have zero projection.

We also use the $\mathcal{C}$-color space $M$, with basis $\{\ee_\alpha: \alpha \text{ is a color of } c\}$. All nonvertical parity equations can then be combined into one equation in $M\otimes W$, where a block with $\mathcal{C}$-color $\alpha$ and signature $w$ contributes $\ee_\alpha\otimes w$.

\begin{lemma}\label{lem:sign}
Let $A,B$ be disjoint two-element column sets and $P,Q$ two-element row sets. Set
$
 i(A,B)=|\{(a,b)\in A\times B:a<b\}|\pmod2.
$
The number of negative-slope edges between $A\times P$ and $B\times Q$ has parity
\begin{equation}\label{eq:signscalar}
 i(A,B)|P\cap Q|.
\end{equation}
More precisely, retaining the $\mathcal{D}$-color and the negative sign but forgetting the $\mathcal{C}$-color gives the parity vector
\begin{equation}\label{eq:signaux}
 i(A,B)\,\pi\ff(P,Q).
\end{equation}
\end{lemma}
\begin{proof}
Let $N$ be the number of pairs $(a,b)\in A\times B$ with $a<b$. For any fixed distinct rows $p\in P$ and $q\in Q$, the number of negative-slope edges is either $N$ or $4-N$, depending on the row order. Both numbers have parity $i(A,B)$. Summing the corresponding $\mathcal{D}$-color vectors proves~\eqref{eq:signaux}. There are $4-|P\cap Q|$ ordered pairs of distinct rows, which gives~\eqref{eq:signscalar}. If $P$ and $Q$ are disjoint, the scalar parity is zero, although the vector in~\eqref{eq:signaux} need not be zero.
\end{proof}

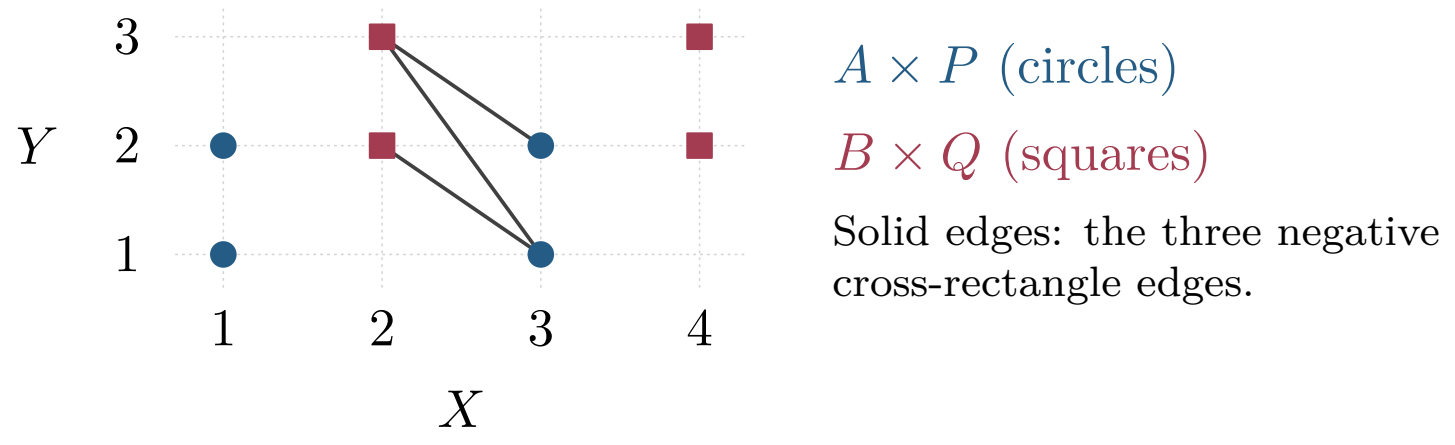
\begin{figure}[!htbp]
\centering
\begin{tikzpicture}[product,x=1.05cm,y=0.72cm]
  \ProdGrid{4}{3}
  \ProdX{1}{1}\ProdX{2}{2}\ProdX{3}{3}\ProdX{4}{4}
  \ProdY{1}{1}\ProdY{2}{2}\ProdY{3}{3}
  \node[anchor=north] at (2.5,-0.05) {$X$};
  \node[anchor=east] at (0.1,2) {$Y$};
  \draw[prod edge] (3,1)--(2,2);
  \draw[prod edge] (3,1)--(2,3);
  \draw[prod edge] (3,2)--(2,3);
  \foreach \xx in {1,3}{\foreach \yy in {1,2}{\ProdPoint{\xx}{\yy}{A}}}
  \foreach \xx in {2,4}{\foreach \yy in {2,3}{\ProdPoint{\xx}{\yy}{B}}}
  \node[anchor=west,text=prodA] at (4.7,2.7) {$A\times P$ (circles)};
  \node[anchor=west,text=prodB] at (4.7,1.9) {$B\times Q$ (squares)};
  \node[anchor=west,font=\scriptsize,align=left] at (4.7,0.95)
    {Solid edges: the three negative\\cross-rectangle edges.};
\end{tikzpicture}
\caption{The displayed coordinate orders are the actual orders on $X$ and $Y$. All eight vertices are shown, and the only graph edges drawn are the three negative-slope cross-rectangle edges.}
\label{fig:eight-vertices}
\end{figure}

\Needspace{10\baselineskip}
For a concrete illustration of~\eqref{eq:signature}--\eqref{eq:signaux}, take the ordered sets
 $X=\{1<2<3<4\}$ and $Y=\{1<2<3\},$
and put
 $$A=\{1,3\}, B=\{2,4\},
 P=\{1,2\},  Q=\{2,3\}\mbox{ and }
 S=(A\times P)\cup(B\times Q).$$
Thus $R_1=R_3=P$ and $R_2=R_4=Q$. Let every column edge have the same $\mathcal{C}$-color $\alpha$, and use three distinct $\mathcal{D}$-colors $r,b,g$ with $d(12)=r$, $d(23)=b$ and $d(13)=g.$ This example illustrates only the parity identities; no uniqueness hypothesis is imposed on $c$ or $d$, see Fig.~\ref{fig:eight-vertices}.

For a block with support $P$ on the left and $Q$ on the right, the four ordered row pairs contribute to~\eqref{eq:signature} as follows:
$$
\begin{array}{c|c|c}
 (p,q)\in P\times Q & \text{edge type and $\mathcal{D}$-color}
                    & \text{contribution to }\ff(P,Q)\\ \hline
 (1,2)&\text{positive diagonal, }r&\ee_{r,+}\\
 (1,3)&\text{positive diagonal, }g&\ee_{g,+}\\
 (2,2)&\text{horizontal}&\hh\\
 (2,3)&\text{positive diagonal, }b&\ee_{b,+}
\end{array}
$$
Consequently, writing $w=\ff(P,Q)$, we obtain
$$
 w=\hh+\ee_{r,+}+\ee_{g,+}+\ee_{b,+}\mbox{ and }
 \ff(Q,P)=\tau w=\hh+\ee_{r,-}+\ee_{g,-}+\ee_{b,-}.
$$
For the two blocks with equal supports, the horizontal terms occur twice and cancel over $\F$. Thus
$$
 \ff(P,P)=2\hh+\ee_{r,+}+\ee_{r,-}=\ttv_r\mbox{ and }
 \ff(Q,Q)=2\hh+\ee_{b,+}+\ee_{b,-}=\ttv_b,
$$
verifying~\eqref{eq:rectanglevector} for both rectangles.

The projection $\pi$ removes the horizontal coordinate and identifies the two signs of each $\mathcal{D}$-color. In this example,
$$
 \pi w=\ee_r+\ee_g+\ee_b\mbox{ and }
 w+\tau w=\ttv_r+\ttv_g+\ttv_b=\ttv_{\pi w}.
$$
Thus~\eqref{eq:taupi} holds: the two copies of $\hh$ cancel, while each $\mathcal{D}$-color retains one positive and one negative coordinate.

To check~\eqref{eq:parityequation}, note that all six column pairs have $\mathcal{C}$-color $\alpha$, with signatures
$$
\begin{array}{c|cccccc}
 \text{column pair}&12&13&14&23&24&34\\ \hline
 \text{signature}&w&\ttv_r&w&\tau w&\ttv_b&w
\end{array}
$$
For instance, block $23$ has support $Q$ on the left and $P$ on the right, so its signature is $\tau w$. Since $3w=w$ over $\F$, the left side of \eqref{eq:parityequation} is
 $$\sum_{x<x'}\ff(R_x,R_{x'})
 =\ttv_r+\ttv_b+3w+\tau w
 =\ttv_r+\ttv_b+(\ttv_r+\ttv_g+\ttv_b)
 =\ttv_g\ne0.$$
Thus~\eqref{eq:parityequation} shows that $S$ is not even-chromatic. Direct counting gives the same conclusion:
$$
\begin{array}{c|cc}
 \text{$\mathcal{D}$-color}&\text{positive diagonal edges}
                         &\text{negative diagonal edges}\\ \hline
 r&4&2\\
 b&4&2\\
 g&3&1
\end{array}
$$
There are also eight horizontal edges of color $(H,\alpha)$ and four vertical edges, with each exact vertical color $(V,\{1,2\})$ and $(V,\{2,3\})$ occurring twice. The total is $8+4+(4+2)+(4+2)+(3+1)=28=\binom82$. Exactly two full colors, $(D,\alpha,g,+)$ and $(D,\alpha,g,-)$, occur an odd number of times, as encoded by $\ttv_g$.

Finally, consider the negative-slope parity in~\eqref{eq:signscalar} and~\eqref{eq:signaux}. Among the four ordered column pairs in $A\times B$, the three pairs $(1,2),(1,4),(3,4)$ satisfy $a<b$, while $(3,2)$ satisfies $a>b$. Hence $i(A,B)=3\bmod2=1$. Also $P\cap Q=\{2\}$. The three nonhorizontal row pairs in $P\times Q$, namely $(1,2),(1,3),(2,3)$, all satisfy $p<q$. Each therefore yields a negative-slope edge exactly for the reversed column pair $(3,2)$. These three edges, drawn above, are
$$
 \{(3,1),(2,2)\}, 
 \{(3,1),(2,3)\}\mbox{ and }
 \{(3,2),(2,3)\}.
$$
Their $\mathcal{D}$-colors are $r,g,b$, respectively. Thus the number $N_-$ of negative-slope edges \emph{between the two rectangles} satisfies
$$
 N_-=3\equiv i(A,B)|P\cap Q|=1\pmod2,
$$
which is \eqref{eq:signscalar}. Retaining their $\mathcal{D}$-colors gives
$
 \ee_r+\ee_g+\ee_b
 =i(A,B)\,\pi\ff(P,Q),
$
as in~\eqref{eq:signaux}. These last two equations concern only cross-rectangle edges. Adding the two internal negative diagonals, of $\mathcal{D}$-colors $r$ and $b$, leaves total negative-slope parity $\ee_g$, in agreement with the full-color calculation above.
\medskip

We now apply these identities to two rectangles. The following lemma will also be used to establish the eight-vertex preservation result needed later.

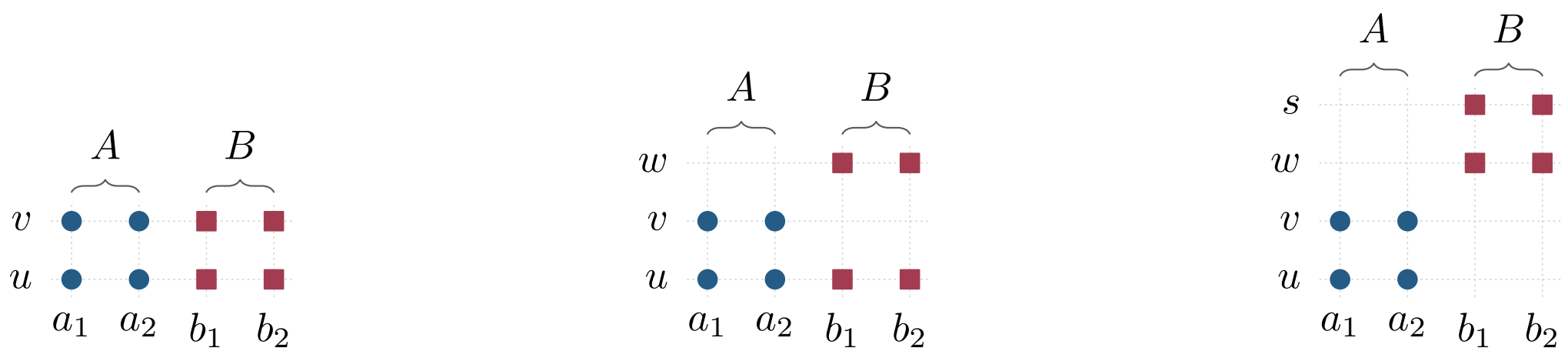
\begin{figure}[!htbp]
\centering
\begin{subfigure}[t]{.32\linewidth}
\centering
\begin{tikzpicture}[product,x=0.58cm,y=0.50cm]
\ProdGrid{4}{2}
\ProdX{1}{a_1}
\ProdX{2}{a_2}
\ProdX{3}{b_1}
\ProdX{4}{b_2}
\ProdY{1}{u}
\ProdY{2}{v}
\ProdPair{1}{2}{2}{A}
\ProdPair{3}{4}{2}{B}

\ProdPoint{1}{1}{A}
\ProdPoint{1}{2}{A}
\ProdPoint{2}{1}{A}
\ProdPoint{2}{2}{A}
\ProdPoint{3}{1}{B}
\ProdPoint{3}{2}{B}
\ProdPoint{4}{1}{B}
\ProdPoint{4}{2}{B}

\end{tikzpicture}
\caption{$P=Q=\{u,v\}$}
\end{subfigure}
\hfill
\begin{subfigure}[t]{.32\linewidth}
\centering
\begin{tikzpicture}[product,x=0.58cm,y=0.50cm]
\ProdGrid{4}{3}
\ProdX{1}{a_1}
\ProdX{2}{a_2}
\ProdX{3}{b_1}
\ProdX{4}{b_2}
\ProdY{1}{u}
\ProdY{2}{v}
\ProdY{3}{w}
\ProdPair{1}{2}{3}{A}
\ProdPair{3}{4}{3}{B}

\ProdPoint{1}{1}{A}
\ProdPoint{1}{2}{A}
\ProdPoint{2}{1}{A}
\ProdPoint{2}{2}{A}
\ProdPoint{3}{1}{B}
\ProdPoint{3}{3}{B}
\ProdPoint{4}{1}{B}
\ProdPoint{4}{3}{B}

\end{tikzpicture}
\caption{$P=\{u,v\},\ Q=\{u,w\}$}
\end{subfigure}
\hfill
\begin{subfigure}[t]{.32\linewidth}
\centering
\begin{tikzpicture}[product,x=0.58cm,y=0.50cm]
\ProdGrid{4}{4}
\ProdX{1}{a_1}
\ProdX{2}{a_2}
\ProdX{3}{b_1}
\ProdX{4}{b_2}
\ProdY{1}{u}
\ProdY{2}{v}
\ProdY{3}{w}
\ProdY{4}{s}
\ProdPair{1}{2}{4}{A}
\ProdPair{3}{4}{4}{B}

\ProdPoint{1}{1}{A}
\ProdPoint{1}{2}{A}
\ProdPoint{2}{1}{A}
\ProdPoint{2}{2}{A}
\ProdPoint{3}{3}{B}
\ProdPoint{3}{4}{B}
\ProdPoint{4}{3}{B}
\ProdPoint{4}{4}{B}

\end{tikzpicture}
\caption{$P=\{u,v\},\ Q=\{w,s\}$}
\end{subfigure}
\caption{The eight vertices in Lemma~\ref{lem:eight}. The three panels show all possible intersection types of the two row supports. Circles belong to $A\times P$ and squares to $B\times Q$.}
\label{fig:two-rectangles}
\end{figure}

\begin{lemma}\label{lem:eight}
If $c,d$ are $K_4$-unique, then $(A\times P)\cup(B\times Q)$ is not even-chromatic for disjoint two-element column sets $A,B$ and arbitrary two-element row sets $P,Q$.
Here Figure~\ref{fig:two-rectangles} displays the three possible row-intersection types.
\end{lemma}

\begin{proof}
Suppose otherwise, and choose a $\mathcal{C}$-color that occurs on exactly one edge of the clique on $A\cup B$. Its edge cannot lie within $A$ or $B$, because the corresponding block signature would be nonzero by~\eqref{eq:rectanglevector}. The edge therefore joins $A$ to $B$, so $\ff(P,Q)=0$. If $|P\cap Q|=1$, the coefficient of $\hh$ is nonzero; if $P=Q$, equation~\eqref{eq:rectanglevector} applies. Hence $P$ and $Q$ are disjoint.

Every block between the two rectangles now has zero signature. Within each rectangle, the vertical and horizontal edges occur in equal-color pairs. After these contributions are discarded, only the two internal diagonal pairs remain. The two positive diagonals must have the same full color, so in particular $d(P)=d(Q)$. On the four rows in $P\cup Q$, the four cross edges have auxiliary parity $\pi\ff(P,Q)=0$, while the two internal edges have the same $\mathcal{D}$-color. This gives an even-chromatic $K_4$ under $d$, contradicting its $K_4$-uniqueness.
\end{proof}

\FloatBarrier

\begin{proposition}\label{prop:eight}
If $c$ is $K_4$-unique and $K_8$-odd, and $d$ is $K_4$-unique, then $c\otimes d$ is $K_8$-odd.
\end{proposition}

\begin{proof}
Suppose that an eight-vertex set is even-chromatic. Its number of odd-sized columns belongs to $\{0,2,4,6,8\}$. By Lemma~\ref{lem:projection}(i), these columns induce an even-chromatic clique under $c$. The values $2$ and $6$ are impossible because $\binom22$ and $\binom62$ are odd. The values $4$ and $8$ are excluded by the $K_4$-uniqueness and $K_8$-oddness of $c$, respectively. Thus all occupied columns have even size.

By Lemma~\ref{lem:projection}(iii), each column has at most four vertices, so the only possible size patterns are
$$
 4+4, 4+2+2\mbox{ and }2+2+2+2.
$$
In the pattern $4+4$, Lemma~\ref{lem:projection}(ii) forces the two supports to coincide. The positive diagonal edges between the columns then have exactly the $\mathcal{D}$-color multiplicities of a $K_4$, contradicting the $K_4$-uniqueness of $d$. In the pattern $4+2+2$, the other two columns cannot repeat all six row pairs in the four-vertex column, again contradicting Lemma~\ref{lem:projection}(ii). Finally, in the pattern $2+2+2+2$, Lemma~\ref{lem:projection}(ii) implies that each two-element row support occurs an even number of times. Pairing equal supports expresses the set as two rectangles, which is excluded by Lemma~\ref{lem:eight}.
\end{proof}

\section{The case for $K_9$}\label{sec:nine}

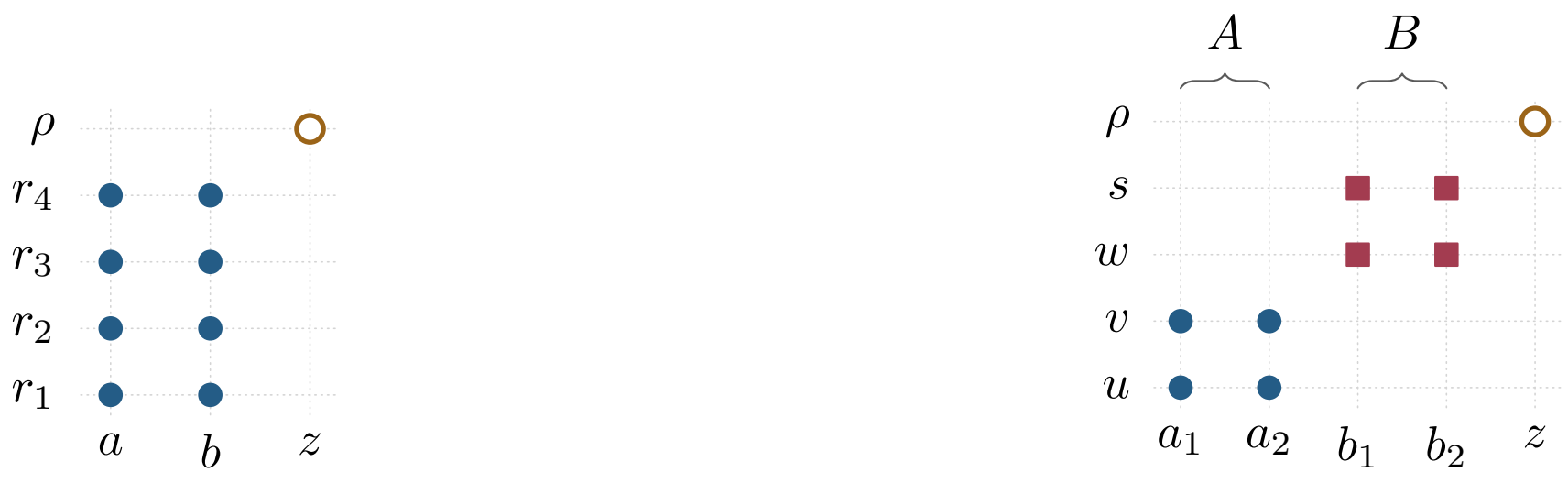
\begin{figure}[!htbp]
\centering
\begin{subfigure}[t]{.46\linewidth}
\centering
\begin{tikzpicture}[product,x=0.72cm,y=0.48cm]
\ProdGrid{3}{5}
\ProdX{1}{a}
\ProdX{2}{b}
\ProdX{3}{z}
\ProdY{1}{r_1}
\ProdY{2}{r_2}
\ProdY{3}{r_3}
\ProdY{4}{r_4}
\ProdY{5}{\rho}

\ProdPoint{1}{1}{A}
\ProdPoint{1}{2}{A}
\ProdPoint{1}{3}{A}
\ProdPoint{1}{4}{A}
\ProdPoint{2}{1}{A}
\ProdPoint{2}{2}{A}
\ProdPoint{2}{3}{A}
\ProdPoint{2}{4}{A}
\ProdPoint{3}{5}{Z}

\end{tikzpicture}
\caption{$(\{a,b\}\times R)\cup\{(z,\rho)\}$}
\end{subfigure}
\hfill
\begin{subfigure}[t]{.50\linewidth}
\centering
\begin{tikzpicture}[product,x=0.64cm,y=0.48cm]
\ProdGrid{5}{5}
\ProdX{1}{a_1}
\ProdX{2}{a_2}
\ProdX{3}{b_1}
\ProdX{4}{b_2}
\ProdX{5}{z}
\ProdY{1}{u}
\ProdY{2}{v}
\ProdY{3}{w}
\ProdY{4}{s}
\ProdY{5}{\rho}
\ProdPair{1}{2}{5}{A}
\ProdPair{3}{4}{5}{B}

\ProdPoint{1}{1}{A}
\ProdPoint{1}{2}{A}
\ProdPoint{2}{1}{A}
\ProdPoint{2}{2}{A}
\ProdPoint{3}{3}{B}
\ProdPoint{3}{4}{B}
\ProdPoint{4}{3}{B}
\ProdPoint{4}{4}{B}
\ProdPoint{5}{5}{Z}

\end{tikzpicture}
\caption{$(A\times P)\cup(B\times Q)\cup\{(z,\rho)\}$}
\end{subfigure}
\caption{The two support types in Lemma~\ref{lem:classification9}, each with nine vertices. Panel (a) uses $R=\{r_1,r_2,r_3,r_4\}$; panel (b) uses $P=\{u,v\}$ and $Q=\{w,s\}$. The hollow vertex is $(z,\rho)$. Distinct row labels are used for clarity: the lemma also permits $\rho$ to coincide with an existing row and, in (b), permits $P,Q$ to intersect or coincide.}
\label{fig:classification9}
\end{figure}

We first classify the column supports of an assumed even-chromatic nine-vertex set and then exclude the two resulting configurations (see Figure~\ref{fig:classification9}). The argument for two rectangles and an additional vertex uses symmetric assumptions on $c$ and $d$, so we may transpose the configuration when its five row coordinates are distinct.

\begin{lemma}\label{lem:classification9}
Suppose that $c$ is $K_5$-odd and $K_9$-odd. If a nine-vertex set $S$ is even-chromatic under $c\otimes d$, then it has one of the following forms:
\begin{enumerate}[label=\textup{(\roman*)},leftmargin=2em]
\item $(\{a,b\}\times R)\cup\{(z,\rho)\}$, where $a,b,z$ are distinct and $|R|=4$;
\item $(A\times P)\cup(B\times Q)\cup\{(z,\rho)\}$, where $A,B$ are disjoint two-element column sets, $z\notin A\cup B$, and $P,Q$ are two-element row sets.
\end{enumerate}
In either form, $\rho$ may coincide with an already occupied row. In \textup{(ii)}, $P$ and $Q$ may also intersect or coincide.
\end{lemma}
\begin{proof}
Since $|S|=9$, the number of odd-sized columns belongs to $\{1,3,5,7,9\}$. By Lemma~\ref{lem:projection}(i), these columns induce an even-chromatic clique under $c$. The values $3$ and $7$ are impossible because $\binom32$ and $\binom72$ are odd, while the hypotheses on $c$ exclude $5$ and $9$. Thus exactly one occupied column has odd size.
By Lemma~\ref{lem:projection}(iii), no column has more than four vertices.

Suppose first that a column has four vertices. There is exactly one odd-sized column, so the possible size patterns are
$$
4+4+1,  4+3+2\mbox{ and } 4+2+2+1.
$$
Each of the six row pairs in the four-vertex column must occur in another column. In the last two patterns, the other columns provide at most $\binom32+\binom22=4$ and $2\binom22=2$ row-pair occurrences, respectively. Neither suffices. In the first pattern, the other four-vertex column must contain all six row pairs, and hence
has the same row support. This gives (i).

If the largest column has size three, there is exactly one three-vertex column, and the size pattern is $3+2+2+2$. Write its support as $\{u,v,w\}$. To repeat its three vertical labels, the other supports must be $\{u,v\}$, $\{u,w\}$ and $\{v,w\}$. Each of these rows then contains three vertices, giving $3\binom32=9$ horizontal edges. This is impossible because the horizontal palette is disjoint from the other palettes and all its color multiplicities must be even.

Finally, if every column has size at most two, there are four two-vertex columns and one singleton. By Lemma~\ref{lem:projection}(ii), each two-element support occurs an even number of times. Pairing columns with equal supports therefore gives \textup{(ii)}.
\end{proof}

To handle the additional vertex, fix its row $\rho$ and write $\pp_U=\ff(\{\rho\},U)$. The two possible column orders give signatures $\pp_U$ and $\tau\pp_U$.

\FloatBarrier

\begin{lemma}\label{lem:star}
Let $|U|=2$.
\begin{enumerate}[label=\textup{(\roman*)},leftmargin=2em]
\item If $\rho\in U$, then $\pp_U$ has coefficient one at $\hh$. If $\rho\notin U$, then its coefficient at $\hh$ is zero.
\item The equality $\pp_U=0$ holds if and only if $\rho\notin U$, $\rho$ lies outside the interval spanned by $U$, and the two edges from $\rho$ to $U$ have the same $\mathcal{D}$-color.
\item If $\rho\notin U$ and $d$ has no monochromatic triangle, then
$ \ttv_U\notin\spn_{\F}\{\pp_U,\tau\pp_U\}.$
\end{enumerate}
\end{lemma}
\begin{proof}
The first two assertions follow directly from the two terms in $\pp_U$. For (iii), write $U=\{u,v\}$ and put
$\gamma=d(\rho u)$, $\delta=d(\rho v)$. We have
$$
\pp_U+\tau\pp_U=\ttv_\gamma+\ttv_\delta.
$$
This vector is zero if $\gamma=\delta$ and otherwise is the sum of two distinct two-sign vectors. In either case, it differs from the single two-sign vector $\ttv_U$. If $\pp_U$ or $\tau\pp_U$ were equal to $\ttv_U$, its two terms would have color $d(uv)$ and opposite signs. In particular, $d(\rho u)=d(\rho v)=d(uv)$, giving a monochromatic triangle. Since $\ttv_U$ is nonzero, none of the four vectors in the indicated span can equal $\ttv_U$.
\end{proof}

We now exclude the first configuration in Lemma~\ref{lem:classification9}.

\begin{lemma}\label{lem:four}
Suppose $c$ is $K_3$-unique and $d$ is $K_4$-odd. Let $a,b,z$ be distinct columns and let $R$ be a four-element row set. For every $\rho\in Y$, the set
$
 (\{a,b\}\times R)\cup\{(z,\rho)\}
$
is not even-chromatic.
\end{lemma}
\begin{proof}
Suppose otherwise, and choose a color that occurs on exactly one edge of the triangle $\{a,b,z\}$ under $c$. If this edge is $ab$, all nonvertical edges with this $\mathcal{C}$-color lie between the two copies of $R$. Each unordered row pair in $R$ determines one positive diagonal and one negative diagonal. Thus, for each $\mathcal{D}$-color $\gamma$, the number of positive diagonals with color coordinates $(c(ab),\gamma)$ equals the multiplicity of $\gamma$ in $K_R$. At least one of these multiplicities is odd because $d$ is $K_4$-odd, a contradiction.

Otherwise, the unique column edge joins $z$ to $a$ or $b$, and its four-edge block has zero signature. The block from $(z,\rho)$ to the other column has the same row supports, so its signature is either the same vector or its $\tau$-image and is also zero. Discarding both blocks leaves the doubled four-set. The positive-diagonal count above gives a color of odd multiplicity in this remaining set, again a contradiction.
\end{proof}

It remains to exclude the second configuration in Lemma~\ref{lem:classification9}. Let $A=\{a_1,a_2\}$ and $B=\{b_1,b_2\}$ be disjoint column sets, let $z\notin A\cup B$, and let $P,Q$ be two-element row sets. For an arbitrary row $\rho$, put
\begin{equation}\label{eq:nineset}
 S=(A\times P)\cup(B\times Q)\cup\{(z,\rho)\}.
\end{equation}
The next lemma provides a cancellation argument used in two of the cases below.

\begin{lemma}\label{lem:cross}
Suppose that $P\cap Q=\emptyset$ and $S$ is even-chromatic. Suppose that a $\mathcal{D}$-color occurs on exactly one edge of $P\cup Q\cup\{\rho\}$, namely $pq$, where $p\in P$, $q\in Q$ and $p,q\ne\rho$. Then every full color occurs an even number of times among the sixteen edges between $A\times P$ and $B\times Q$. Moreover, every $\mathcal{C}$-color occurs an even number of times among the four edges between $A$ and $B$.
\end{lemma}

\begin{proof}
Let $\delta=d(pq)$. Its only occurrences as a diagonal $\mathcal{D}$-color coordinate on $S$ are the four edges $(a,p)(b,q)$, where $a\in A$ and $b\in B$. The condition $p,q\ne\rho$ excludes edges incident to $(z,\rho)$. For each $\mathcal{C}$-color $\alpha$, let $N_\alpha^+$ count the pairs $(a,b)\in A\times B$ with $c(ab)=\alpha$ and $a<b$, and let $N_\alpha^-$ count those with $c(ab)=\alpha$ and $a>b$. If $p<q$, these numbers are the multiplicities of $(D,\alpha,\delta,+)$ and $(D,\alpha,\delta,-)$, respectively; if $p>q$, their roles are reversed.

Even-chromaticity therefore implies that $N_\alpha^+$ and $N_\alpha^-$ are both even. The contribution of the cross-rectangle edges to the $\alpha$-summand is consequently
$$
N_\alpha^+\ff(P,Q)+N_\alpha^-\tau\ff(P,Q)=0.
$$
This proves the first assertion, regardless of any coincidences among the other $\mathcal{D}$-colors. The second follows because $N_\alpha^++N_\alpha^-$ is even for every $\alpha$.
\end{proof}

\begin{figure}[!htbp]
\centering
\begin{subfigure}[t]{.48\linewidth}
\centering
\begin{tikzpicture}[product,x=0.64cm,y=0.49cm]
\ProdGrid{5}{5}
\ProdX{1}{a_1}
\ProdX{2}{a_2}
\ProdX{3}{b_1}
\ProdX{4}{b_2}
\ProdX{5}{z}
\ProdY{1}{p}
\ProdY{2}{p'}
\ProdY{3}{q}
\ProdY{4}{q'}
\ProdY{5}{\rho}
\ProdPair{1}{2}{5}{A}
\ProdPair{3}{4}{5}{B}
\draw[prod focus] (1,1) to[bend left=14] (3,3);
\draw[prod focus] (1,1) -- (4,3);
\draw[prod focus] (2,1) -- (3,3);
\draw[prod focus] (2,1) -- (4,3);
\ProdPoint{1}{1}{A}
\ProdPoint{1}{2}{A}
\ProdPoint{2}{1}{A}
\ProdPoint{2}{2}{A}
\ProdPoint{3}{3}{B}
\ProdPoint{3}{4}{B}
\ProdPoint{4}{3}{B}
\ProdPoint{4}{4}{B}
\ProdPoint{5}{5}{Z}

\end{tikzpicture}
\caption{$\rho\notin P\cup Q$}
\end{subfigure}
\hfill
\begin{subfigure}[t]{.48\linewidth}
\centering
\begin{tikzpicture}[product,x=0.64cm,y=0.49cm]
\ProdGrid{5}{4}
\ProdX{1}{a_1}
\ProdX{2}{a_2}
\ProdX{3}{b_1}
\ProdX{4}{b_2}
\ProdX{5}{z}
\ProdY{1}{p}
\ProdY{2}{p'=\rho}
\ProdY{3}{q}
\ProdY{4}{q'}
\ProdPair{1}{2}{4}{A}
\ProdPair{3}{4}{4}{B}
\draw[prod focus] (1,1) to[bend left=14] (3,3);
\draw[prod focus] (1,1) -- (4,3);
\draw[prod focus] (2,1) -- (3,3);
\draw[prod focus] (2,1) -- (4,3);
\ProdPoint{1}{1}{A}
\ProdPoint{1}{2}{A}
\ProdPoint{2}{1}{A}
\ProdPoint{2}{2}{A}
\ProdPoint{3}{3}{B}
\ProdPoint{3}{4}{B}
\ProdPoint{4}{3}{B}
\ProdPoint{4}{4}{B}
\ProdPoint{5}{2}{Z}

\end{tikzpicture}
\caption{$\rho\in P\setminus\{p\}$}
\end{subfigure}
\caption{The cancellation mechanism in Lemma~\ref{lem:cross}, with $P=\{p,p'\}$, $Q=\{q,q'\}$ and $\delta=d(pq)$ unique on the occupied row set. The four highlighted edges are the only diagonal occurrences of color $\delta$. Panel (a) illustrates the case where $\rho\notin P\cup Q$. Panel (b) illustrates an allowed row coincidence; the case $\rho\in Q\setminus\{q\}$ is symmetric. The argument cancels all sixteen cross-rectangle edges, although only these four diagnostic edges are drawn.}
\label{fig:cancellation}
\end{figure}
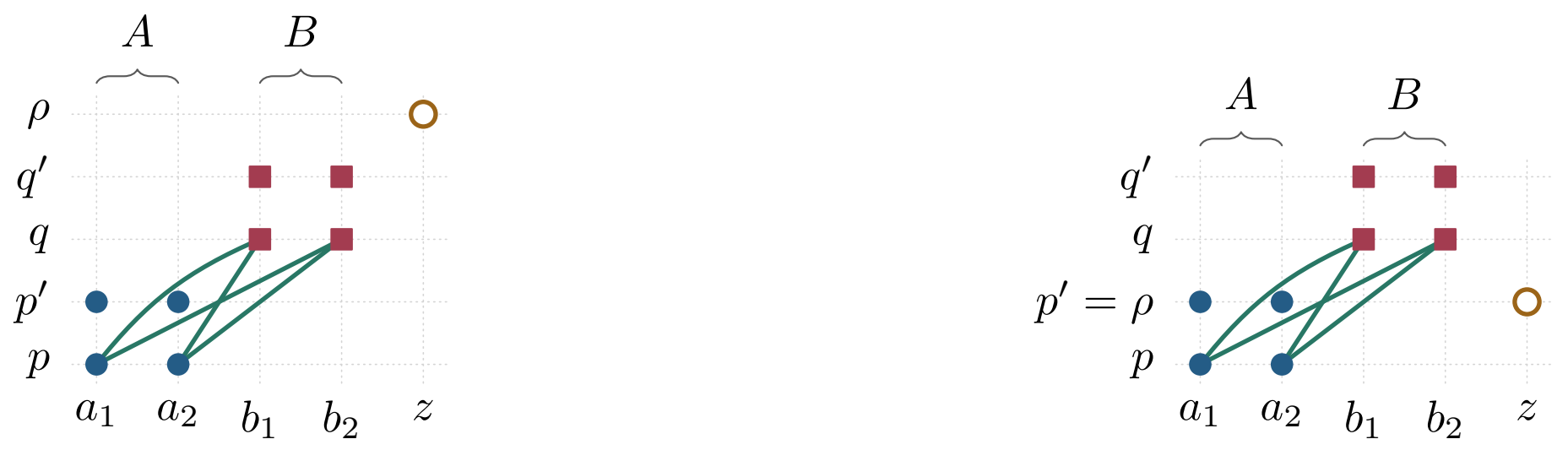

Figure~\ref{fig:cancellation} isolates the four edges used to prove the cancellation.

\begin{lemma}\label{lem:transpose}
Suppose $P,Q,\{\rho\}$ are pairwise disjoint. If the set $S$ in~\eqref{eq:nineset} is even-chromatic under $c\otimes d$, then its transpose
$
S^{\mathsf T}=(P\times A)\cup(Q\times B)\cup\{(\rho,z)\}
$
is even-chromatic under $d\otimes c$.
\end{lemma}
\begin{proof}
Transposition sends a diagonal color $(D,\alpha,\gamma,\varepsilon)$ to $(D,\gamma,\alpha,\varepsilon)$. This is a bijective relabelling, so diagonal color parities are preserved.

The four horizontal edges of $S$ form two equal-color pairs, one on each exact column pair $A$ and $B$. After transposition, these become two pairs of vertical edges with exact labels $A$ and $B$. Likewise, the four vertical edges form two pairs with exact row labels $P$ and $Q$; after transposition, each pair becomes a pair of horizontal edges with the same color. Thus all horizontal and vertical color multiplicities are even both before and after transposition.
\end{proof}

\begin{proposition}\label{prop:nine}
If $c$ and $d$ both have property $\U$, then no set of the form~\eqref{eq:nineset} is even-chromatic.
\end{proposition}
\begin{proof}
Suppose that $S$ is even-chromatic. Every vertical color already occurs evenly: each row pair $P$ or $Q$ occurs twice, or four times if $P=Q$. We therefore consider only the nonvertical parity equations.

Since $c$ is $K_5$-unique, some color occurs on exactly one edge among the five columns. This edge cannot be internal to $A$ or $B$, since its signature would be the nonzero vector $\ttv_P$ or $\ttv_Q$. If it joins $z$ to a column in $A$, Observation~\ref{obs:uniquezero} gives $\pp_P=0$, and hence $\rho\notin P$. Both blocks from $(z,\rho)$ to $A\times P$ then have zero signature, regardless of their $\mathcal{C}$-colors or column orders. The same argument applies to $B$ and $\pp_Q$.

\medskip\noindent\textbf{Case 1: $P=Q$.}
Every edge among the four doubled columns has signature $\ttv_P\ne0$. Thus the unique $\mathcal{C}$-color on the five columns must lie on an edge incident to $z$, forcing $\pp_P=0$. We may discard all eight edges incident to $(z,\rho)$. The remaining union of two rectangles would then be even-chromatic, contradicting Lemma~\ref{lem:eight}.

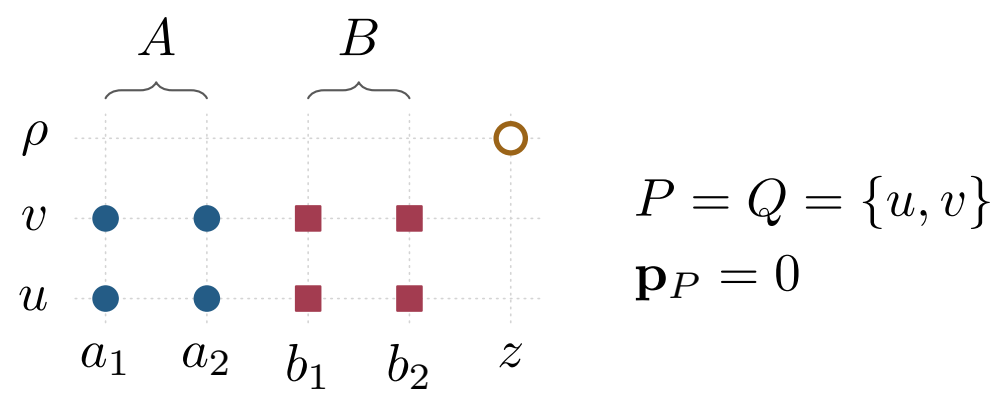
\begin{figure}[!htbp]
\centering
\begin{tikzpicture}[product]
  \ProdGrid{5}{3}
  \ProdX{1}{a_1}\ProdX{2}{a_2}\ProdX{3}{b_1}\ProdX{4}{b_2}\ProdX{5}{z}
  \ProdY{1}{u}\ProdY{2}{v}\ProdY{3}{\rho}
  \ProdPair{1}{2}{3}{A}\ProdPair{3}{4}{3}{B}
  \foreach \xx in {1,2}{\foreach \yy in {1,2}{\ProdPoint{\xx}{\yy}{A}}}
  \foreach \xx in {3,4}{\foreach \yy in {1,2}{\ProdPoint{\xx}{\yy}{B}}}
  \ProdPoint{5}{3}{Z}
  \node[anchor=west,align=left] at (6,2.2) {$P=Q=\{u,v\}$};
  \node[anchor=west,align=left] at (6,1.25) {$\pp_P=0$};
\end{tikzpicture}
\caption{Proposition~\ref{prop:nine}, Case 1. The four doubled columns have the same row support. Under the assumed evenness, the unique $\mathcal{C}$-color forces $\pp_P=0$; in particular $\rho\notin P$. All eight edges incident to $(z,\rho)$ can then be discarded.}
\label{fig:case1}
\end{figure}
\FloatBarrier

\medskip\noindent\textbf{Case 2: $|P\cap Q|=1$.}
Write $P=\{u,v\}$ and $Q=\{u,w\}$.
The coefficient of $\hh$ in $\ff(P,Q)$ is one, so the unique $\mathcal{C}$-color cannot lie on an edge between $A$ and $B$. It must lie on an edge incident to $z$. Exchanging the rectangles if necessary, we may assume that
\begin{equation}\label{eq:case2star}
 \pp_P=0,  \rho\notin P\mbox{ and } d(\rho u)=d(\rho v).
\end{equation}
Discard the two zero-signature blocks from $(z,\rho)$ to $A\times P$.

Now choose a color $c(xx')=\alpha$ that occurs on exactly one edge among the four columns $A\cup B$, and write this edge with $x<x'$. Its block signature $\ff(R_x,R_{x'})$ must be canceled by the remaining blocks from $(z,\rho)$ to $B\times Q$. In particular, this signature belongs to
\begin{equation}\label{eq:case2span}
 \spn_{\F}\{\pp_Q,\tau\pp_Q\}.
\end{equation}
By~\eqref{eq:case2star}, either $\rho\notin P\cup Q$ or $\rho=w$.

Figure~\ref{fig:case2} displays these two possibilities.

\begin{figure}[!htbp]
\centering
\begin{subfigure}[t]{.47\linewidth}
\centering
\begin{tikzpicture}[product]
  \ProdGrid{5}{4}
  \ProdX{1}{a_1}\ProdX{2}{a_2}\ProdX{3}{b_1}\ProdX{4}{b_2}\ProdX{5}{z}
  \ProdY{1}{u}\ProdY{2}{v}\ProdY{3}{w}\ProdY{4}{\rho}
  \ProdPair{1}{2}{4}{A}\ProdPair{3}{4}{4}{B}
  \foreach \xx in {1,2}{\foreach \yy in {1,2}{\ProdPoint{\xx}{\yy}{A}}}
  \foreach \xx in {3,4}{\foreach \yy in {1,3}{\ProdPoint{\xx}{\yy}{B}}}
  \ProdPoint{5}{4}{Z}
\end{tikzpicture}
\caption{Case 2a: $\rho\notin P\cup Q$.}
\end{subfigure}\hfill
\begin{subfigure}[t]{.47\linewidth}
\centering
\begin{tikzpicture}[product]
  \path[use as bounding box] (-.3,-.15) rectangle (5.35,5.03);
  \ProdGrid{5}{3}
  \ProdX{1}{a_1}\ProdX{2}{a_2}\ProdX{3}{b_1}\ProdX{4}{b_2}\ProdX{5}{z}
  \ProdY{1}{u}\ProdY{2}{v}\ProdY{3}{w=\rho}
  \ProdPair{1}{2}{3}{A}\ProdPair{3}{4}{3}{B}
  \draw[prod focus] (3,3) to[bend left=20] (5,3);
  \draw[prod focus] (4,3)--(5,3);
  \foreach \xx in {1,2}{\foreach \yy in {1,2}{\ProdPoint{\xx}{\yy}{A}}}
  \foreach \xx in {3,4}{\foreach \yy in {1,3}{\ProdPoint{\xx}{\yy}{B}}}
  \ProdPoint{5}{3}{Z}
\end{tikzpicture}
\caption{Case 2b: $\rho=w$.}
\end{subfigure}
\caption{Proposition~\ref{prop:nine}, Case 2: $P=\{u,v\}$ and $Q=\{u,w\}$ share one row. In both panels $\pp_P=0$, so the singleton-to-$A$ blocks are discarded. The singleton-to-$B$ signatures have no horizontal term in (a), whereas each has a horizontal term in (b). The two horizontal edges from the singleton are highlighted in (b).}
\label{fig:case2}
\end{figure}
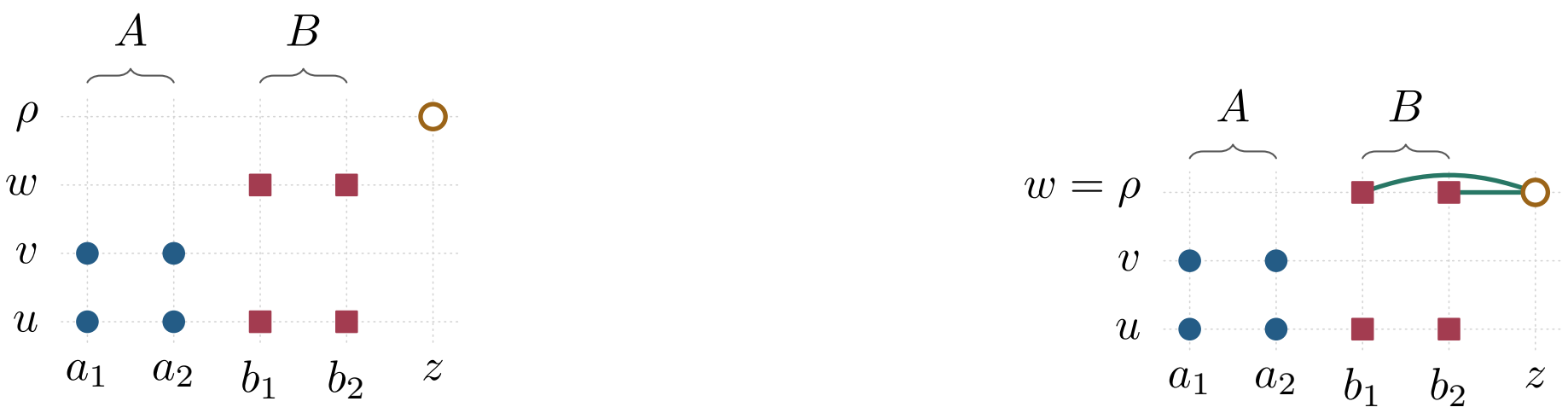
\FloatBarrier

\smallskip\noindent\textbf{Case 2a: $\rho\notin P\cup Q$.}
Every vector in~\eqref{eq:case2span} has zero $\hh$-coordinate. Thus $xx'$ cannot join $A$ to $B$. Lemma~\ref{lem:star}(iii) excludes the internal $B$-block, whose signature is $\ttv_Q$. Suppose instead that the signature is $\ttv_P$. The sum $\pp_Q+\tau\pp_Q$ is either zero or a sum of two distinct two-sign vectors, and therefore cannot equal $\ttv_P$. If $\pp_Q$ or $\tau\pp_Q$ equals $\ttv_P$, then $d(\rho u)=d(\rho w)=d(uv)$. Together with~\eqref{eq:case2star}, this gives $d(\rho u)=d(\rho v)=d(uv)$, a monochromatic triangle. All possible locations of $xx'$ are thus excluded.

\smallskip\noindent\textbf{Case 2b: $\rho=w$.}
Put $\kappa=d(wu)=d(wv)$ and $\lambda=d(uv)$. Since $d$ has no monochromatic triangle, $\lambda\ne\kappa$. The vector $\pp_Q$ consists of $\hh$ and one term of color $\kappa$. Every vector in~\eqref{eq:case2span} is therefore supported on $\hh,\ee_{\kappa,+},\ee_{\kappa,-}$. On the other hand, $\ttv_P$ has color $\lambda$, and
$$
\ff(P,Q)=\hh
 +\ee_{\kappa,\sgn(w-u)}
 +\ee_{\lambda,\sgn(u-v)}
 +\ee_{\kappa,\sgn(w-v)}
$$
has a nonzero $\lambda$-coordinate. Neither $\ff(P,Q)$ nor its $\tau$-image can belong to~\eqref{eq:case2span}. The unique column edge must therefore be $b_1b_2$, with signature $\ttv_Q$. Each remaining singleton block has coefficient one at $\hh$, so neither can cancel $\ttv_Q$ alone. Both must belong to the $\mathcal{C}$-color summand $c(b_1b_2)$, giving
$
c(zb_1)=c(zb_2)=c(b_1b_2).
$
This contradicts the absence of monochromatic triangles in $c$.

\medskip\noindent\textbf{Case 3: $P\cap Q=\emptyset$ and
$\rho\in P\cup Q$.}
After exchanging the rectangles if necessary, write $P=\{u,v\}$, $\rho=u$ and $Q=\{w,s\}$. We first establish a fact that does not depend on the location of the unique $\mathcal{C}$-color.

\begin{claim}\label{claim:internalrows}
No color that occurs on exactly one edge of $\{u,v,w,s\}$ under $d$ can lie on the internal row pair $uv$ or $ws$.
\end{claim}
\begin{proof}
If the unique row pair is $ws$, its only nonvertical occurrences as a $\mathcal{D}$-color coordinate are the two diagonals of $B\times Q$. These diagonals have opposite signs, so each has a unique full color, contradicting even-chromaticity.

If the unique row pair is $uv$, put $\delta=d(uv)$ and $\alpha=c(a_1a_2)$. Its only nonvertical occurrences are the two diagonals of $A\times P$ and the edges from $(z,u)$ to $(a_1,v)$ and $(a_2,v)$. The diagonals have opposite signs, so the two additional edges must match them, one for each sign. In particular, $c(za_1)=c(za_2)=\alpha$, giving a monochromatic triangle in $c$.
\end{proof}

Choose a color that occurs on exactly one edge among the five columns. Its edge cannot join $z$ to $A$, since this would produce a unique $\mathcal{C}$-color. It also cannot join $A$ to $B$: such an edge would force $\ff(P,Q)=0$. A unique $\mathcal{D}$-color on the four rows could then not lie on a cross row edge, since it would contribute exactly one term to $\ff(P,Q)$. It would have to lie on an internal pair, contradicting Claim~\ref{claim:internalrows}. Thus the unique column edge joins $z$ to $B$, and
\begin{equation}\label{eq:case3star}
 \pp_Q=0\mbox{ and } d(uw)=d(us).
\end{equation}
Discard both blocks from $(z,u)$ to $B\times Q$.

Now choose a unique $\mathcal{D}$-color on the four rows. By~\eqref{eq:case3star}, it does not lie on $uw$ or $us$, and Claim~\ref{claim:internalrows} excludes the internal pairs. It must therefore lie on $vw$ or $vs$. Neither endvertex is $\rho=u$, so Lemma~\ref{lem:cross} allows us to discard all sixteen cross-rectangle edges.

Put $\alpha=c(a_1a_2)$ and $\beta=c(b_1b_2)$. By Lemma~\ref{lem:cross}, every $\mathcal{C}$-color occurs evenly among the four cross column edges. If $\alpha=\beta$, the entire four-column clique would be even-chromatic, contradicting $K_4$-uniqueness. Thus $\alpha\ne\beta$. The remaining parity equations contain $\ttv_P$ in the $\alpha$-summand, $\ttv_Q$ in the $\beta$-summand, and the two blocks from $(z,u)$ to $A\times P$. Each of the last two blocks has coefficient one at $\hh$. Since no other remaining term has an $\hh$-coordinate, these blocks must have the same $\mathcal{C}$-color. Their contributions lie in a single summand and therefore cannot cancel the two nonzero vectors in the distinct $\alpha$- and $\beta$-summands. This is a contradiction.

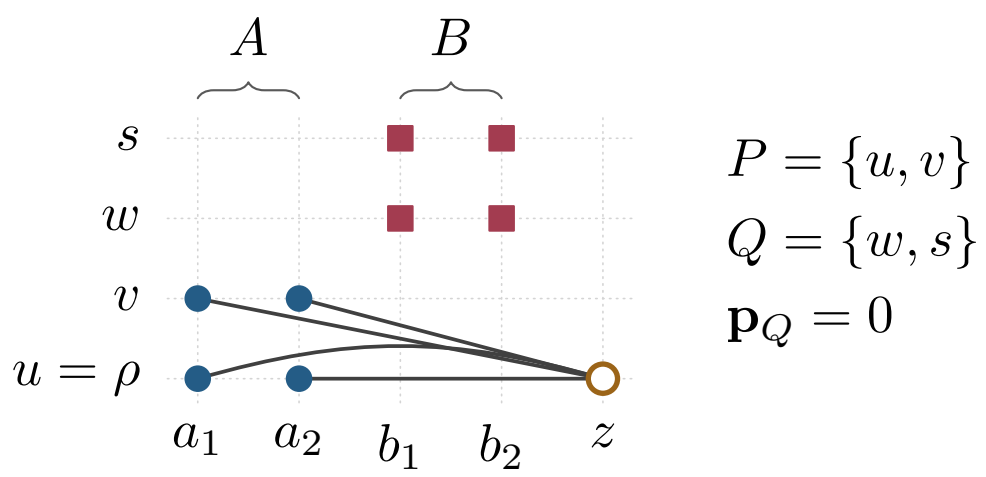
\begin{figure}[!htbp]
\centering
\begin{tikzpicture}[product]
  \ProdGrid{5}{4}
  \ProdX{1}{a_1}\ProdX{2}{a_2}\ProdX{3}{b_1}\ProdX{4}{b_2}\ProdX{5}{z}
  \ProdY{1}{u=\rho}\ProdY{2}{v}\ProdY{3}{w}\ProdY{4}{s}
  \ProdPair{1}{2}{4}{A}\ProdPair{3}{4}{4}{B}
  \draw[prod edge] (5,1) to[bend right=16] (1,1);
  \draw[prod edge] (5,1)--(2,1);
  \draw[prod edge] (5,1)--(1,2);
  \draw[prod edge] (5,1)--(2,2);
    \foreach \xx in {1,2}{\foreach \yy in {1,2}{\ProdPoint{\xx}{\yy}{A}}}
  \foreach \xx in {3,4}{\foreach \yy in {3,4}{\ProdPoint{\xx}{\yy}{B}}}
  \ProdPoint{5}{1}{Z}
  \node[anchor=west,align=left] at (6,3.2) {$P=\{u,v\}$\\[3pt]$Q=\{w,s\}$};
  \node[anchor=west,align=left] at (6,1.7) {$\pp_Q=0$};
 
\end{tikzpicture}
\caption{Proposition~\ref{prop:nine}, Case 3: disjoint row supports with $\rho=u\in P$. The displayed edges are the two remaining singleton-to-$A$ blocks, each with one horizontal term. The proof first discards the singleton-to-$B$ blocks and then all edges between the two rectangles. The internal rectangle contributions $\ttv_P$ and $\ttv_Q$ remain in distinct $\mathcal{C}$-color summands.}
\label{fig:case3}
\end{figure}
\FloatBarrier

\medskip\noindent\textbf{Case 4: $P,Q,\{\rho\}$ are pairwise disjoint.}
There are five distinct columns and five distinct rows. Choose a unique color in each of these five-vertex cliques. The chosen column edge cannot be internal to $A$ or $B$, since the two diagonals in that rectangle would have distinct, unique full colors. The same argument excludes an internal row edge in $P$ or $Q$. Thus each chosen edge either joins the two pairs or is incident to the singleton coordinate.

The two chosen edges cannot both join the pairs: a unique column edge between $A$ and $B$ and a unique row edge between $P$ and $Q$ would together select exactly one diagonal edge of $S$, whose full color would be unique. At least one chosen edge is therefore incident to the singleton coordinate. If necessary, apply Lemma~\ref{lem:transpose} to transpose the grid and exchange $c$ and $d$. Both colorings satisfy $\U$, so the assumptions remain valid. After also exchanging the rectangles if necessary, we may assume that the chosen unique column edge joins $z$ to $A$. Hence $\pp_P=0$, and both blocks from $(z,\rho)$ to $A\times P$ can be discarded.

The chosen unique row edge cannot join $\rho$ to $P$, since together with the chosen unique column edge it would again determine a single diagonal edge with a unique full color. This leaves two possibilities.

\begin{figure}[!htbp]
\centering
\begin{tikzpicture}[product]
  \ProdGrid{5}{5}
  \ProdX{1}{a_1}\ProdX{2}{a_2}\ProdX{3}{b_1}\ProdX{4}{b_2}\ProdX{5}{z}
  \ProdY{1}{u}\ProdY{2}{v}\ProdY{3}{w}\ProdY{4}{s}\ProdY{5}{\rho}
  \ProdPair{1}{2}{5}{A}\ProdPair{3}{4}{5}{B}
  \foreach \xx in {1,2}{\foreach \yy in {1,2}{\ProdPoint{\xx}{\yy}{A}}}
  \foreach \xx in {3,4}{\foreach \yy in {3,4}{\ProdPoint{\xx}{\yy}{B}}}
    \ProdPoint{5}{5}{Z}
  \node[anchor=west,align=left] at (6,3.6) {$P=\{u,v\}$\\[3pt]$Q=\{w,s\}$};
  \node[anchor=west,align=left] at (6,1.8) {$P,Q,\{\rho\}$ are pairwise disjoint};
\end{tikzpicture}
\caption{Proposition~\ref{prop:nine}, Case 4: nine vertices on five distinct columns and five distinct rows. Transposition exchanges the roles of the two coordinates. After transposing and exchanging the rectangles if needed, a unique $\mathcal{C}$-color joins $z$ to $A$, forcing $\pp_P=0$. The two possible locations of the unique row color are shown in Figure~\ref{fig:case4ab}.}
\label{fig:case4}
\end{figure}
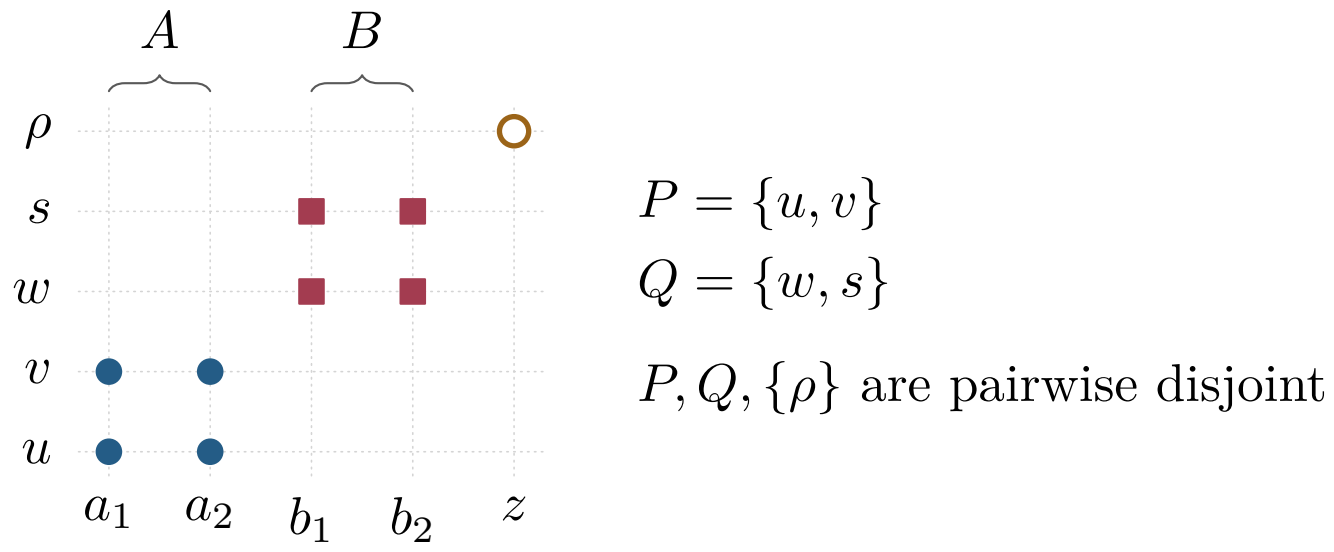
\FloatBarrier

\smallskip\noindent\textbf{Case 4a: the unique row edge joins $\rho$ to $Q$.}
Write this edge as $\rho w$. Its $\mathcal{D}$-color appears on exactly two edges of $S$, namely $(z,\rho)(b_1,w)$ and $(z,\rho)(b_2,w)$. Their full colors must agree, so
$$
c(zb_1)=c(zb_2),
$$
and their signs also agree. Equality of signs implies that $z$ lies outside the column interval spanned by $b_1$ and $b_2$. Consequently, for the other row $s\in Q$, the edges $(z,\rho)(b_1,s)$ and $(z,\rho)(b_2,s)$ also have the same $\mathcal{C}$-color, $\mathcal{D}$-color and sign. Every full color therefore occurs evenly among the four edges from $(z,\rho)$ to $B\times Q$. We have now discarded all singleton edges, contradicting Lemma~\ref{lem:eight}.

\smallskip\noindent\textbf{Case 4b: the unique row edge joins $P$ to $Q$.}
Its endvertices differ from $\rho$, so Lemma~\ref{lem:cross} allows all cross-rectangle edges to be discarded. As in Case 3, the internal $\mathcal{C}$-colors $\alpha=c(a_1a_2)$ and $\beta=c(b_1b_2)$ are distinct, since otherwise the four-column clique would be even-chromatic. The remaining terms are $\ttv_P$ in the $\alpha$-summand, $\ttv_Q$ in the $\beta$-summand, and the two singleton-to-$B$ blocks. Each nonzero internal vector requires a contribution in its own summand. Since there are only two remaining blocks and each belongs to a single summand, there must be one block in each of the $\alpha$- and $\beta$-summands. In particular, exactly one block must cancel $\ttv_Q$. Its signature is $\pp_Q$ or $\tau\pp_Q$, contrary to Lemma~\ref{lem:star}(iii).

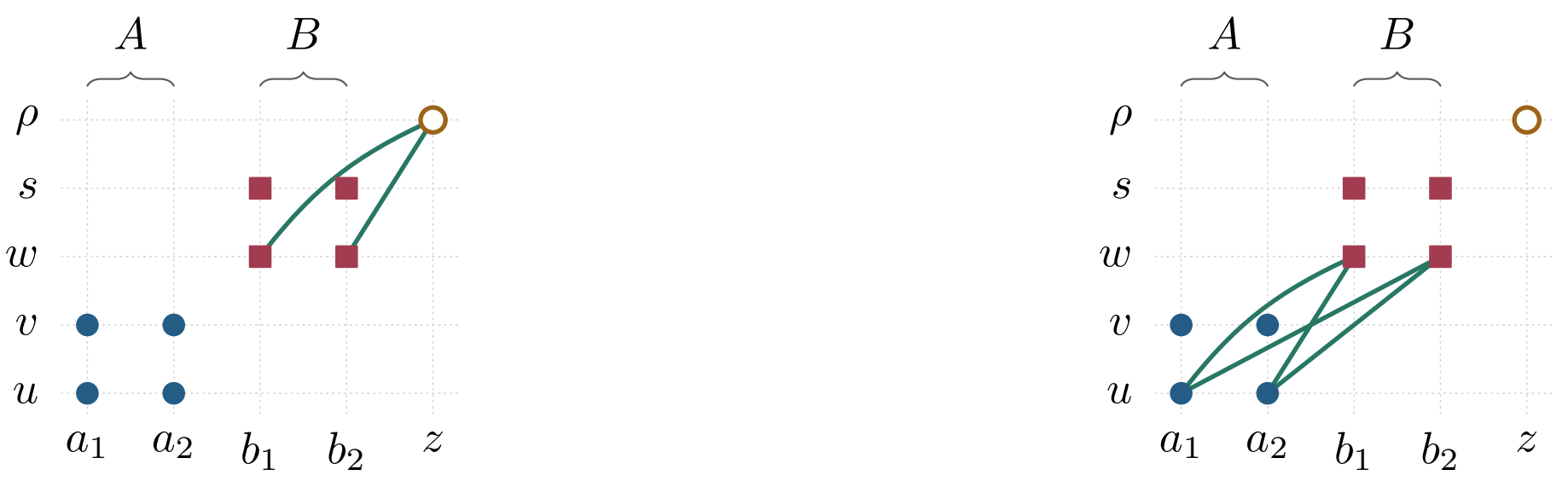
\begin{figure}[!htbp]
\centering
\begin{subfigure}[t]{.47\linewidth}
\centering
\begin{tikzpicture}[product]
  \ProdGrid{5}{5}
  \ProdX{1}{a_1}\ProdX{2}{a_2}\ProdX{3}{b_1}\ProdX{4}{b_2}\ProdX{5}{z}
  \ProdY{1}{u}\ProdY{2}{v}\ProdY{3}{w}\ProdY{4}{s}\ProdY{5}{\rho}
  \ProdPair{1}{2}{5}{A}\ProdPair{3}{4}{5}{B}
  \draw[prod focus] (5,5) to[bend right=14] (3,3);
  \draw[prod focus] (5,5)--(4,3);
  \foreach \xx in {1,2}{\foreach \yy in {1,2}{\ProdPoint{\xx}{\yy}{A}}}
  \foreach \xx in {3,4}{\foreach \yy in {3,4}{\ProdPoint{\xx}{\yy}{B}}}
  \ProdPoint{5}{5}{Z}
\end{tikzpicture}
\caption{Case 4a: the unique row pair is $\rho w$.}
\end{subfigure}\hfill
\begin{subfigure}[t]{.47\linewidth}
\centering
\begin{tikzpicture}[product]
  \ProdGrid{5}{5}
  \ProdX{1}{a_1}\ProdX{2}{a_2}\ProdX{3}{b_1}\ProdX{4}{b_2}\ProdX{5}{z}
  \ProdY{1}{u}\ProdY{2}{v}\ProdY{3}{w}\ProdY{4}{s}\ProdY{5}{\rho}
  \ProdPair{1}{2}{5}{A}\ProdPair{3}{4}{5}{B}
  \draw[prod focus] (1,1) to[bend left=14] (3,3);
  \draw[prod focus] (1,1)--(4,3);
  \draw[prod focus] (2,1)--(3,3);
  \draw[prod focus] (2,1)--(4,3);
  \foreach \xx in {1,2}{\foreach \yy in {1,2}{\ProdPoint{\xx}{\yy}{A}}}
  \foreach \xx in {3,4}{\foreach \yy in {3,4}{\ProdPoint{\xx}{\yy}{B}}}
  \ProdPoint{5}{5}{Z}
\end{tikzpicture}
\caption{Case 4b: the unique row pair is $uw$.}
\end{subfigure}
\caption{The two subcases of Case 4 in Proposition~\ref{prop:nine}. Highlighted edges are exactly the diagonal occurrences of the chosen unique $\mathcal{D}$-color: two singleton-to-$B$ edges in (a), and four cross-rectangle edges in (b). Their parity equations allow all singleton-to-$B$ edges to be discarded in (a), and all cross-rectangle edges to be discarded in (b).}
\label{fig:case4ab}
\end{figure}
\FloatBarrier

These four cases exhaust all possibilities for $P,Q,\rho$.
\end{proof}

Combining the support classification with these exclusions gives the required preservation statement.

\FloatBarrier

\begin{proposition}\label{prop:closure9}
If $c$ has property $\U$ and is $K_9$-odd, and $d$ has
property $\U$, then $c\otimes d$ has property $\U$ and is
$K_9$-odd.
\end{proposition}
\begin{proof}
Proposition~\ref{prop:unique} preserves $\U$. Suppose that $S$ is an even-chromatic nine-vertex set. Since $K_5$-uniqueness implies $K_5$-oddness, Lemma~\ref{lem:classification9} applies. Lemma~\ref{lem:four} excludes its first configuration, using $K_3$-uniqueness of $c$ and $K_4$-oddness of $d$. Proposition~\ref{prop:nine} excludes its second configuration. Hence no such $S$ exists.
\end{proof}

\begin{proof}[\bf Proof of Theorem~\ref{thm:k9}]
Put
$
\mathcal P_9=\{K_3\uniq,K_4\uniq,K_5\uniq,K_9\odd\}.
$
These conditions all hold on one vertex and are preserved by restriction. By Proposition~\ref{prop:closure9}, $c\otimes d$ has property $\mathcal P_9$ whenever $c$ has property $\mathcal P_9$ and $d$ has property $\U$. Lemma~\ref{lem:iteration} therefore gives a $K_9$-odd coloring of $K_n$ with at most $\exp(C_9(\log n)^{2/3})$ colors, proving Theorem~\ref{thm:k9}. This coloring is also $K_s$-unique for $3\le s\le5$.
\end{proof}

\section{The case for $K_{12}$}\label{sec:twelve}

We first classify the column supports of an assumed even-chromatic twelve-vertex set. Two configurations remain: a doubled four-set together with a rectangle, and three rectangles. A count of positive diagonal edges excludes the first. To exclude the second, we analyze two possible patterns of row incidence. We then prove the preservation of $K_6$-uniqueness needed to complete the induction.

\begin{lemma}\label{lem:classification12}
Suppose that $c$ is $K_4$-odd, $K_8$-odd and $K_{12}$-odd.
If a twelve-vertex set $S$ is even-chromatic under $c\otimes d$, then it has one of the following forms:
\begin{enumerate}[label=\textup{(\roman*)},leftmargin=2em]
\item $(A\times R)\cup(B\times P)$, where $A,B$ are disjoint two-element column sets, $|R|=4$ and $|P|=2$;
\item $\bigcup_{i=1}^3(A_i\times P_i)$, where the $A_i$ are pairwise disjoint two-element column sets and $|P_i|=2$ for each $i$.
\end{enumerate}
The row sets may intersect.
\end{lemma}

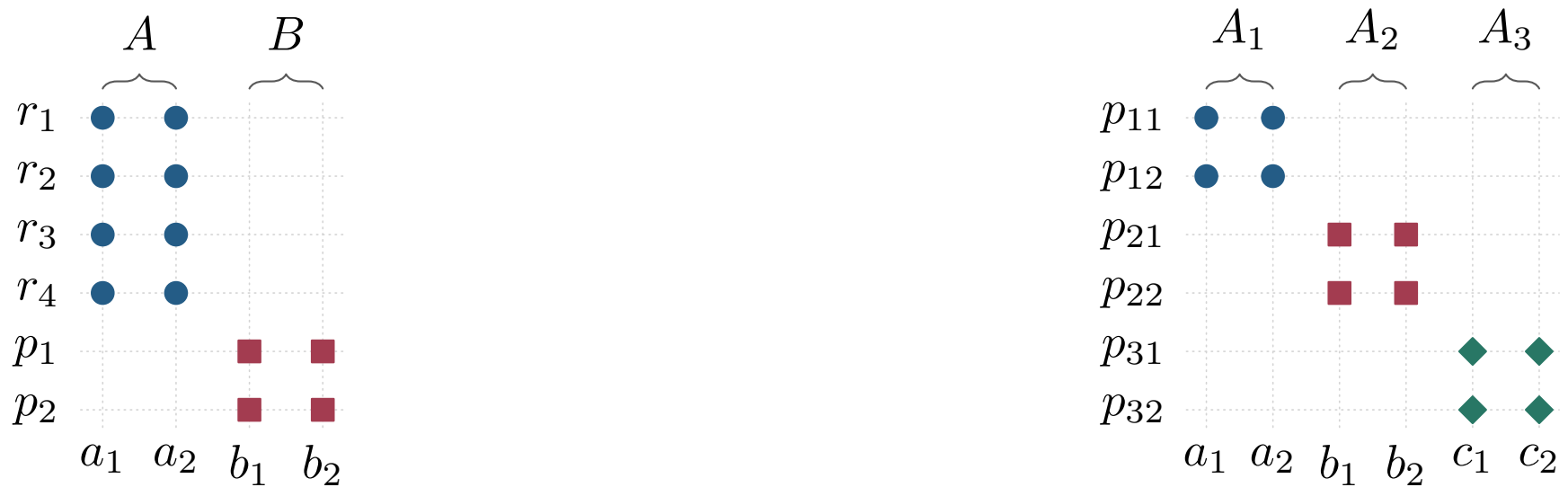
\begin{figure}[!htbp]
  \centering
  \begin{subfigure}[t]{0.47\textwidth}
    \centering
    \begin{tikzpicture}[product,x=0.54cm,y=0.43cm]
      \ProdGrid{4}{6}
      \ProdX{1}{a_1}\ProdX{2}{a_2}
      \ProdX{3}{b_1}\ProdX{4}{b_2}
      \foreach \yy/\lab in {6/r_1,5/r_2,4/r_3,3/r_4,2/p_1,1/p_2}
        {\ProdY{\yy}{\lab}}
      \ProdPair{1}{2}{6}{A}\ProdPair{3}{4}{6}{B}
      \foreach \xx in {1,2}{\foreach \yy in {3,4,5,6}{\ProdPoint{\xx}{\yy}{A}}}
      \foreach \xx in {3,4}{\foreach \yy in {1,2}{\ProdPoint{\xx}{\yy}{B}}}
    \end{tikzpicture}
    \caption{$(A\times R)\cup(B\times P)$,\\ $|R|=4$, $|P|=2$.}
  \end{subfigure}\hfill
  \begin{subfigure}[t]{0.47\textwidth}
    \centering
    \begin{tikzpicture}[product,x=0.49cm,y=0.43cm]
      \ProdGrid{6}{6}
      \ProdX{1}{a_1}\ProdX{2}{a_2}
      \ProdX{3}{b_1}\ProdX{4}{b_2}
      \ProdX{5}{c_1}\ProdX{6}{c_2}
      \foreach \yy/\lab in {6/p_{11},5/p_{12},4/p_{21},3/p_{22},2/p_{31},1/p_{32}}
        {\ProdY{\yy}{\lab}}
      \ProdPair{1}{2}{6}{A_1}\ProdPair{3}{4}{6}{A_2}\ProdPair{5}{6}{6}{A_3}
      \foreach \xx in {1,2}{\foreach \yy in {5,6}{\ProdPoint{\xx}{\yy}{A}}}
      \foreach \xx in {3,4}{\foreach \yy in {3,4}{\ProdPoint{\xx}{\yy}{B}}}
      \foreach \xx in {5,6}{\foreach \yy in {1,2}{\ProdPoint{\xx}{\yy}{C}}}
    \end{tikzpicture}
    \caption{$\bigcup_{i=1}^3(A_i\times P_i)$,\\ $P_i=\{p_{i1},p_{i2}\}$.}
  \end{subfigure}
  \caption{The two column-support configurations in Lemma~\ref{lem:classification12}. Each panel contains twelve occupied product coordinates. The row supports are drawn on distinct levels for clarity; the lemma allows different supports to share row coordinates, including coincident two-element supports.}
  \label{fig:classification12}
\end{figure}

Figure~\ref{fig:classification12} illustrates these two configurations.

\begin{proof}
The number of odd-sized columns belongs to
$\{0,2,4,6,8,10,12\}$. By Lemma~\ref{lem:projection}(i), these columns induce an even-chromatic clique under $c$. The values $2,6,10$ are impossible because the corresponding numbers of edges are odd; the values $4,8,12$ are excluded by the hypotheses on $c$. Thus every occupied column has positive even size. By Lemma~\ref{lem:projection}(iii), their sizes are at most six.

If a column has size six, the possible column-size distributions are
$$
6+6, 6+4+2\mbox{ and } 6+2+2+2.
$$
Every row pair occurs in an even number of column supports by
Lemma~\ref{lem:projection}(ii). The six-vertex column contains fifteen row pairs. In the last two distributions, the other columns supply at most seven and three row-pair occurrences, respectively, so neither is possible. In the first distribution, the two supports must coincide. The two columns then determine exactly $\binom62=15$ positive diagonal edges, whose full-color multiplicities cannot all be even.

If the largest column has size four, the possible distributions are
$$
4+4+4,  4+4+2+2\mbox{ and } 4+2+2+2+2.
$$
The last is impossible because the columns outside the four-vertex column supply only four row-pair occurrences. To exclude $4+4+4$, write the three supports as $R_1,R_2,R_3$. Their clique edge sets satisfy
$$
E(K_{R_1})\mathbin\triangle E(K_{R_2})=E(K_{R_3}).
$$
Taking cardinalities gives $12-2\binom{|R_1\cap R_2|}{2}=6$, so
$|R_1\cap R_2|=3$. The symmetric difference is therefore a
$K_{2,3}$ on five vertices, whereas the right-hand side is a $K_4$, a contradiction.

For the distribution $4+4+2+2$, let $R_1,R_2$ be the four-element supports. If they differ, the symmetric difference of their clique edge sets has at least six edges. It cannot equal the symmetric difference of the two row-pair edges supplied by the remaining columns. Hence $R_1=R_2$, and the two two-element supports are equal as well. This gives (i).

Finally, if every column has size two, there are six occupied columns. By Lemma~\ref{lem:projection}(ii), every two-element support occurs an even number of times. Pairing equal supports gives (ii).
\end{proof}

We now turn to three rectangles. The next two lemmas address the row-incidence patterns that require us to retain individual colors. The column and row orders remain arbitrary throughout.
\phantomsection\label{sec:three}

\FloatBarrier

\begin{lemma}\label{lem:pathfour}
Let $A_1,A_2,A_3$ be pairwise disjoint two-element column sets, and let
$
P_1=\{u,v\},  P_2=\{w,s\}\mbox{ and } P_3=\{u,w\},
$
where $u,v,w,s$ are distinct. Suppose that $c$ and $d$ are $K_4$-unique and $\ff(P_1,P_2)=0$. Then $\bigcup_{i=1}^3(A_i\times P_i)$ is not
even-chromatic.
\end{lemma}

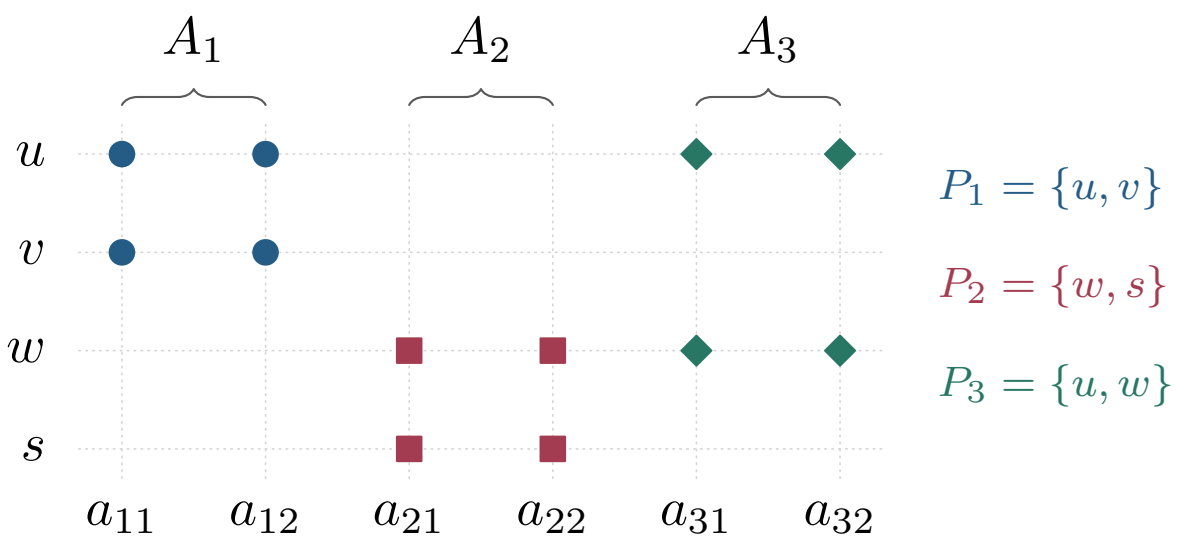
\begin{figure}[!htbp]
  \centering
  \begin{tikzpicture}[product,x=0.95cm,y=0.65cm]
    \ProdGrid{6}{4}
    \ProdX{1}{a_{11}}\ProdX{2}{a_{12}}
    \ProdX{3}{a_{21}}\ProdX{4}{a_{22}}
    \ProdX{5}{a_{31}}\ProdX{6}{a_{32}}
    \ProdY{4}{u}\ProdY{3}{v}\ProdY{2}{w}\ProdY{1}{s}
    \ProdPair{1}{2}{4}{A_1}\ProdPair{3}{4}{4}{A_2}\ProdPair{5}{6}{4}{A_3}
    \foreach \xx in {1,2}{\foreach \yy in {3,4}{\ProdPoint{\xx}{\yy}{A}}}
    \foreach \xx in {3,4}{\foreach \yy in {1,2}{\ProdPoint{\xx}{\yy}{B}}}
    \foreach \xx in {5,6}{\foreach \yy in {2,4}{\ProdPoint{\xx}{\yy}{C}}}
    \node[prod note,anchor=west,text=prodA] at (6.6,3.65) {$P_1=\{u,v\}$};
    \node[prod note,anchor=west,text=prodB] at (6.6,2.65) {$P_2=\{w,s\}$};
    \node[prod note,anchor=west,text=prodC] at (6.6,1.65) {$P_3=\{u,w\}$};
  \end{tikzpicture}
  \caption{The twelve-vertex configuration in Lemma~\ref{lem:pathfour}, with $A_i=\{a_{i1},a_{i2}\}$. Shared row coordinates are aligned: rectangles $1$ and $3$ meet row $u$, while rectangles $2$ and $3$ meet row $w$.}
  \label{fig:pathfour}
\end{figure}

The shared rows in this configuration are aligned in Figure~\ref{fig:pathfour}.

\begin{proof}
Suppose that the set is even-chromatic. Write $\alpha_i=c(A_i)$, and define vectors in the $\mathcal{C}$-color space $M$ by
$$
C_{ij}=\sum_{a\in A_i,b\in A_j}\ee_{c(ab)}\mbox{ and } 
 O_{ij}=\sum_{\substack{a\in A_i,b\in A_j\\a<b}}\ee_{c(ab)}.
$$
If $C_{ij}=0$, reversing the column orientation gives the same parity vector $O_{ij}$. Thus each fixed pair of distinct rows, one from each support, contributes $O_{ij}$ to the positive-slope $\mathcal{C}$-color parity, regardless of the order of those rows.

Choose a $\mathcal{D}$-color unique on the four rows. Its edge cannot join $P_1$ to $P_2$, since $\ff(P_1,P_2)=0$. After interchanging the first two rectangles and relabelling the rows if necessary, we may assume that $\delta=d(uv)$ is unique on these four rows. The only diagonal edges of $\mathcal{D}$-color $\delta$ are the two internal diagonals of the first rectangle and the four edges joining its $v$-row to the $u$-row of the third rectangle. Summing the full-color parity equations over both signs gives $C_{13}=0$. Since the first rectangle has one positive and one negative diagonal, the equation for either sign then gives
 $O_{13}=\ee_{\alpha_1}.$

Only the first and third rectangles share row $u$, and only the second and third share row $w$. The two internal horizontal edges of each rectangle cancel. The remaining horizontal parity equation is therefore $C_{13}+C_{23}=0$, giving $C_{23}=0$.

Now retain only the $\mathcal{C}$-color and the positive-slope tag. The three rectangles contribute $\ee_{\alpha_1}+\ee_{\alpha_2}+\ee_{\alpha_3}$. Every block between the first and second rectangles vanishes by the assumed signature equation. Each of the other two rectangle pairs has exactly three ordered pairs of distinct rows. Since $C_{13}=C_{23}=0$, these pairs contribute $O_{13}$ and $O_{23}$, respectively.
We obtain
$$
\begin{aligned}
0=\ee_{\alpha_1}+\ee_{\alpha_2}+\ee_{\alpha_3}+O_{13}+O_{23} \mbox{ and }
O_{23}=\ee_{\alpha_2}+\ee_{\alpha_3}.
\end{aligned}
$$
If $\alpha_2=\alpha_3$, the equation $C_{23}=0$ makes the $K_4$ on $A_2\cup A_3$ even-chromatic. Otherwise, the equation for $O_{23}$ shows that both colors occur among its cross edges, and $C_{23}=0$ makes each cross multiplicity positive and even. There are only four cross edges, so each color occurs twice. Including the internal edges, this $K_4$ has two colors, each occurring three times. In either case, there is no uniquely colored edge, contradicting the $K_4$-uniqueness of $c$.
\end{proof}
\FloatBarrier

\Needspace{10\baselineskip}
\begin{lemma}\label{lem:paththree}
Suppose that $c$ is both $K_4$-unique and $K_6$-unique, and that $d$ is $K_4$-unique. Let $A,B,C$ be pairwise disjoint two-element column sets, and let $p_1,p_2,p_3,q_1,q_2$ be distinct rows. Put
$
P=\{p_1,p_2\},  Q=\{p_2,p_3\}\mbox{ and } 
 R=\{q_1,q_2\}.
$
Then the twelve-vertex set
$
(A\times P)\cup(B\times Q)\cup(C\times R)
$
is not even-chromatic under $c\otimes d$.
\end{lemma}

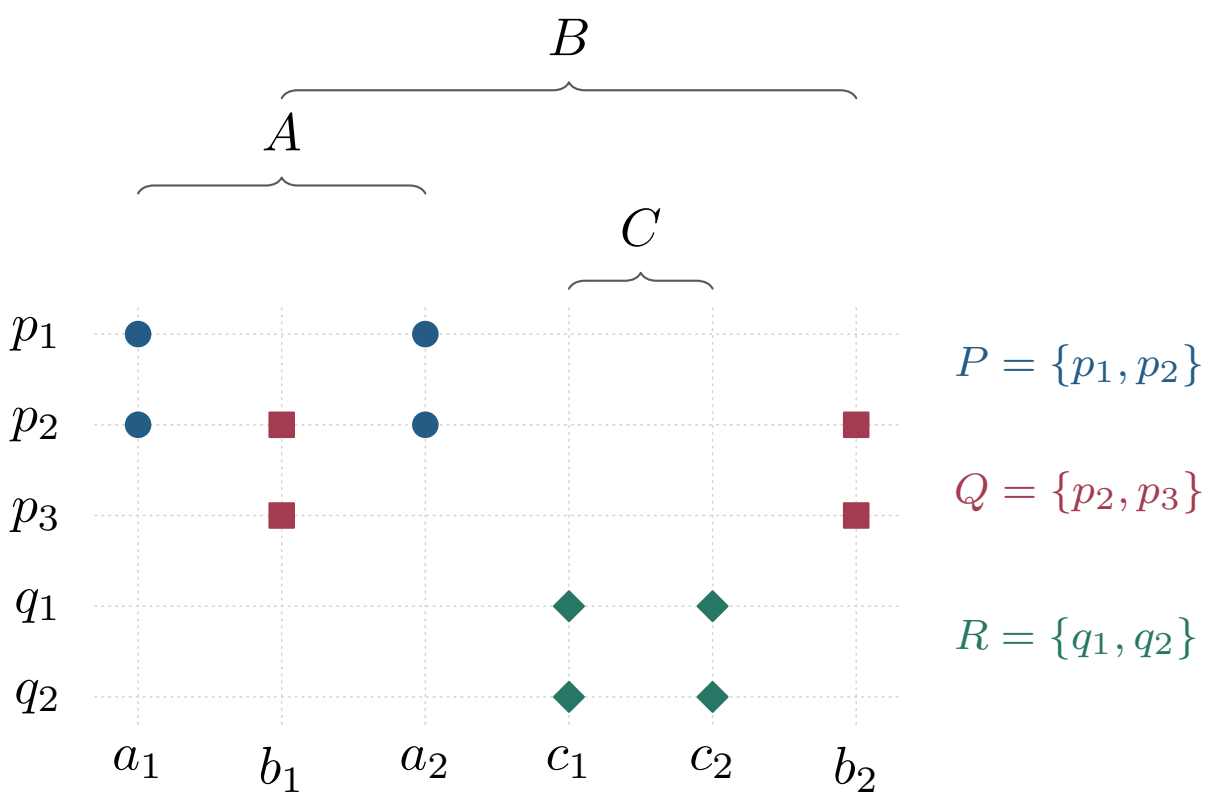
\begin{figure}[!htbp]
  \centering
  \begin{tikzpicture}[product,x=0.95cm,y=0.60cm]
    \ProdGrid{6}{5}
    \ProdX{1}{a_1}\ProdX{2}{b_1}\ProdX{3}{a_2}
    \ProdX{4}{c_1}\ProdX{5}{c_2}\ProdX{6}{b_2}
    \ProdY{5}{p_1}\ProdY{4}{p_2}\ProdY{3}{p_3}
    \ProdY{2}{q_1}\ProdY{1}{q_2}
    \ProdPair{4}{5}{5}{C}
    \ProdPair{1}{3}{6.05}{A}
    \ProdPair{2}{6}{7.1}{B}
    \foreach \xx in {1,3}{\foreach \yy in {4,5}{\ProdPoint{\xx}{\yy}{A}}}
    \foreach \xx in {2,6}{\foreach \yy in {3,4}{\ProdPoint{\xx}{\yy}{B}}}
    \foreach \xx in {4,5}{\foreach \yy in {1,2}{\ProdPoint{\xx}{\yy}{C}}}
    \node[prod note,anchor=west,text=prodA] at (6.6,4.65) {$P=\{p_1,p_2\}$};
    \node[prod note,anchor=west,text=prodB] at (6.6,3.25) {$Q=\{p_2,p_3\}$};
    \node[prod note,anchor=west,text=prodC] at (6.6,1.65) {$R=\{q_1,q_2\}$};
  \end{tikzpicture}
  \caption{The configuration in Lemma~\ref{lem:paththree}: $(A\times P)\cup(B\times Q)\cup(C\times R)$, with twelve occupied product coordinates, where $A=\{a_1,a_2\}$, $B=\{b_1,b_2\}$ and $C=\{c_1,c_2\}$. Only $P$ and $Q$ share a row, namely $p_2$; $R$ is disjoint from both. The row pairs form the path $p_1-p_2-p_3$ and the disjoint edge $q_1q_2$. Under the assumed evenness, the proof obtains $\ff(P,R)=0$ after a possible exchange of $A$ and $B$. The illustrative column order $a_1<b_1<a_2<c_1<c_2<b_2$ has $A,B$ interlaced and $C$ inside $B$, consistent with $i(A,B)=1$ and $i(B,C)=0$ derived in the proof; this order is not assumed in the lemma.}
  \label{fig:paththree}
\end{figure}

Figure~\ref{fig:paththree} illustrates the product configuration and its shared row.

\begin{proof}
Suppose that this set is even-chromatic. Choose a $\mathcal{C}$-color unique on the six columns. Its edge cannot be internal to $A$, $B$ or $C$, by~\eqref{eq:rectanglevector}. Nor can it join $A$ to $B$, because $\ff(P,Q)$ has coefficient one at $\hh$. It therefore joins $C$ to $A$ or $B$. Exchanging $A,B$ and $p_1,p_3$ if necessary, we may assume that it joins $A$ to $C$. Its full-color parity equation gives
\begin{equation}\label{eq:p3edge-zero-block}
 \ff(P,R)=0.
\end{equation}
In particular, all four blocks between $A$ and $C$ have zero signature, irrespective of their $\mathcal{C}$-colors and orientations.

Using the $\mathcal{D}$-color space and the projection $\pi$ defined in Section~\ref{sec:parity}, write
$$
a=d(p_1p_2),  b=d(p_2p_3), 
 z=d(p_1p_3),  x=d(q_1q_2) \mbox{ and }
 \Delta=\pi\ff(Q,R).
$$
Since $Q$ and $R$ are disjoint, $\Delta$ is the parity vector of four $\mathcal{D}$-colors and therefore has even Hamming weight.

For two disjoint column pairs $G,H$, recall
$$
i(G,H)=\bigl|\{(g,h)\in G\times H:g<h\}\bigr|\pmod2.
$$
This number is one exactly when the two pairs interlace. Indeed, the two possible interlacing orders have one and three inversions, whereas each of the other four orders has zero, two or four. We first show that
\begin{equation}\label{eq:p3edge-order-color}
 i(A,B)=1,  i(B,C)=0,  z=x \mbox{ and } \Delta\ne0.
\end{equation}

Consider the $\mathcal{D}$-color parity vector of all negative-slope edges, forgetting their $\mathcal{C}$-colors. Each rectangle contributes its internal $\mathcal{D}$-color once. For a pair of rectangles on column sets $G,H$, Lemma~\ref{lem:sign} gives $i(G,H)$ times the projected signature of their respective row supports. The $AC$ contribution vanishes by
\eqref{eq:p3edge-zero-block}. The vanishing of the total vector therefore gives
\begin{equation}\label{eq:p3edge-negative}
 \ee_a+\ee_b+\ee_x
 +i(A,B)(\ee_a+\ee_b+\ee_z)
 +i(B,C)\Delta=0.
\end{equation}
Summing all coordinates of this equation yields $1+i(A,B)=0$, because $\Delta$ has even weight. Thus $i(A,B)=1$, and
\eqref{eq:p3edge-negative} reduces to
\begin{equation}\label{eq:p3edge-negative-reduced}
 \ee_z+\ee_x+i(B,C)\Delta=0.
\end{equation}

The projected signatures $\pi\ff$ are additive in their row sets over $\F$. Since $P\triangle Q=\{p_1,p_3\}$,
\eqref{eq:p3edge-zero-block} implies
$$
\pi\ff(\{p_1,p_3\},R)
 =\pi\ff(P,R)+\pi\ff(Q,R)=\Delta.
$$
Hence the total color parity vector of $d$ on the four rows
$\{p_1,p_3,q_1,q_2\}$ is
$
\ee_z+\ee_x+\Delta.
$
It is nonzero because $d$ is $K_4$-unique. If $i(B,C)=1$, equation~\eqref{eq:p3edge-negative-reduced} would make this vector zero. Hence $i(B,C)=0$, and the same equation gives $z=x$. The displayed nonzero vector is consequently $\Delta$. This proves
\eqref{eq:p3edge-order-color}.

Similarly, by~\eqref{eq:p3edge-zero-block}, the $\mathcal{D}$-color parity vector on $P\cup R$ is $\ee_a+\ee_x$. Since $d$ is $K_4$-unique, this vector is nonzero, and therefore
\begin{equation}\label{eq:p3edge-a-x}
 a\ne x.
\end{equation}

We now retain the $\mathcal{C}$-colors, using the $\mathcal{C}$-color space $M$. For an ordered pair of column sets $G,H$, define
$$
X_{GH}^{+}
 =\sum_{\substack{g\in G,\ h\in H\\g<h}}\ee_{c(gh)}\mbox{ and }
 X_{GH}^{-}
 =\sum_{\substack{g\in G,\ h\in H\\g>h}}\ee_{c(gh)}.
$$
The superscripts specify column order, rather than the slope of an individual grid edge. Write $\alpha,\beta,\gamma$ for the colors of the internal column edges of $A,B,C$, respectively.

Only the $AB$ blocks have a horizontal term. Their horizontal parity equations give $X_{AB}^{+}=X_{AB}^{-}$. Since $i(A,B)=1$, one orientation class contains exactly one column edge. Thus, for some $\mathcal{C}$-color $\kappa$,
\begin{equation}\label{eq:p3edge-X}
 X_{AB}^{+}=X_{AB}^{-}=\ee_\kappa.
\end{equation}
Apply $\pi$ in each $\mathcal{C}$-color coordinate to the full parity equation. The internal rectangle signatures map to zero, the $AC$ signatures vanish, and the two $AB$ contributions cancel by
\eqref{eq:p3edge-X}. Consequently
$
(X_{BC}^{+}+X_{BC}^{-})\otimes\Delta=0.
$
Since $\Delta\ne0$, there is a vector $Y\in M$ such that
\begin{equation}\label{eq:p3edge-Y}
 X_{BC}^{+}=X_{BC}^{-}=Y.
\end{equation}
In particular, each $\mathcal{C}$-color occurs evenly among the four cross edges of $AB$, and likewise among those of $BC$. If $\alpha=\beta$, the $K_4$ on $A\cup B$ would be even-chromatic under $c$; if $\beta=\gamma$, the same would hold on $B\cup C$. Hence $\beta\neq \alpha, \gamma$.

Recall that the linear map $v\mapsto\ttv_v$ is injective. We have
$$
 \ff(P,Q)+\tau\ff(P,Q)=\ttv_a+\ttv_b+\ttv_x \mbox{ and }
 \ff(Q,R)+\tau\ff(Q,R)=\ttv_\Delta,
$$
where we have used $z=x$. By
\eqref{eq:p3edge-zero-block}, \eqref{eq:p3edge-X} and
\eqref{eq:p3edge-Y}, the full diagonal parity equations are therefore equivalent to the single tensor equation
\begin{equation}\label{eq:p3edge-tensor}
 (\ee_\alpha+\ee_\kappa)\otimes\ttv_a
 +(\ee_\beta+\ee_\kappa)\otimes\ttv_b +(\ee_\gamma+\ee_\kappa)\otimes\ttv_x
 +Y\otimes\ttv_\Delta=0.
\end{equation}

Suppose first that $\kappa\ne\beta$. By $\beta\neq \alpha, \gamma$, the $\beta$ coordinate of
\eqref{eq:p3edge-tensor} is
 $\ttv_b=Y_\beta\ttv_\Delta.$
This is impossible: the left side corresponds to a single $\mathcal{D}$-color, whereas $\Delta$ is nonzero and has even weight. Thus $\kappa=\beta$. The same coordinate now gives
$
 \ttv_a+\ttv_x=Y_\beta\ttv_\Delta.
$
By \eqref{eq:p3edge-a-x}, it follows that
\begin{equation}\label{eq:p3edge-Delta}
 Y_\beta=1\mbox{ and }  \Delta=\ee_a+\ee_x.
\end{equation}
If $\alpha\ne\gamma$, the $\alpha$ coordinate of
\eqref{eq:p3edge-tensor} would give
$
 \ttv_a=Y_\alpha(\ttv_a+\ttv_x),
$
again impossible. Hence $\alpha=\gamma$. The $\alpha$ coordinate, together with \eqref{eq:p3edge-Delta}, gives $Y_\alpha=1$. All remaining $\mathcal{C}$-color coordinates give zero, so
 $Y=\ee_\alpha+\ee_\beta.$

Each orientation class of the four $BC$ column edges must therefore contain at least two edges, one of color $\alpha$ and one of color $\beta$. There are only four such edges altogether, so each class consists of exactly these two edges. The internal edges of $B,C$ have colors $\beta,\alpha$, respectively. Thus the $K_4$ on $B\cup C$ has precisely two colors, each occurring three times. This contradicts the $K_4$-uniqueness of $c$ and completes the proof.
\end{proof}

\FloatBarrier

\begin{proposition}\label{prop:three}
Suppose $c$ is $K_4$-unique and $K_6$-unique, and $d$ is $K_4$-unique. Let $A_1,A_2,A_3$ be pairwise disjoint two-element column sets and let $P_1,P_2,P_3$ be two-element row sets. Then $S=\bigcup_{i=1}^3(A_i\times P_i)$ is not even-chromatic.
\end{proposition}
\begin{proof}
Suppose otherwise and choose a $\mathcal{C}$-color unique on the six columns.
Its edge cannot be internal to an $A_i$, by
\eqref{eq:rectanglevector}. If it joins $A_i$ to $A_j$, it forces
$\ff(P_i,P_j)=0$. The horizontal coordinate excludes an intersection
of size one, and~\eqref{eq:rectanglevector} excludes equality.
Thus $P_i,P_j$ are disjoint.

Each rectangle contributes exactly one negative diagonal, and hence the three rectangles contribute three negative diagonals. By Lemma~\ref{lem:sign}, if all pairwise intersections of the row pairs have even size, then the number of negative diagonals between distinct rectangles is even. Therefore, the total number of negative-slope edges is odd. This contradicts evenness of the negative palette. Thus some two row supports meet in exactly one vertex. View the three row pairs as edges on the set of rows. They contain both a disjoint pair and a pair meeting in one vertex. The three row pairs are distinct: if two coincided, the third could not form both a disjoint pair and a pair meeting in exactly one row. A simple graph with three edges and these two properties is either a path of length three or a path of length two together with a disjoint edge.

In the first case, the disjoint pair is unique, so the uniquely colored column edge forces the zero signature required in Lemma~\ref{lem:pathfour}. Lemma~\ref{lem:paththree} excludes the second case. Both give a contradiction.
\end{proof}

We can now combine the support classification with these exclusion arguments to prove local preservation of $K_{12}$-oddness.

\begin{proposition}\label{prop:twelve}
Suppose $c$ is $K_4$-unique, $K_6$-unique, $K_8$-odd and
$K_{12}$-odd, and $d$ is $K_4$-unique. Then $c\otimes d$ is
$K_{12}$-odd.
\end{proposition}
\begin{proof}
Suppose that a twelve-vertex set is even-chromatic. By Lemma~\ref{lem:classification12}, it has one of the two forms listed there.

Consider $(A\times R)\cup(B\times P)$ with $|R|=4$, $|P|=2$. Choose a color unique under $c$ on $A\cup B$. If its edge is internal to $A$, the positive diagonal part of its block has the $\mathcal{D}$-color multiplicities of $d$ on $R$. These are not all even, by $K_4$-uniqueness of $d$. If its edge is internal to $B$, its block has signature $\ttv_{d(P)}\ne0$. Therefore the edge must join $A$ to $B$, forcing $\ff(R,P)=0$. All four blocks between the two parts now vanish, regardless of their column orders. The two remaining internal blocks contain $\binom42+\binom22=7$ positive diagonals in total, an odd number, so their full-color multiplicities cannot all be even. Proposition~\ref{prop:three} excludes the second form.
\end{proof}

To complete the induction for Theorem~\ref{thm:k12}, we also need to preserve $K_6$-uniqueness.

\begin{proposition}\label{prop:unique6}
Suppose that $c$ is $K_3$-unique, $K_4$-unique and $K_6$-unique, and that $d$ is $K_3$-unique and $K_4$-unique. Then $c\otimes d$ is $K_6$-unique.
\end{proposition}

\begin{proof}
Suppose that a set $S$ of six grid vertices has no uniquely occurring color. We first restrict its possible column sizes. If a column contains the rows $u$ and $v$, then the vertical color indexed by the exact pair $\{u,v\}$ must occur in a second column. In particular, if a column contains $q\geq 2$ vertices, each of its $q$ row coordinates must occur at least once outside that column. Hence $2q\leq 6$, so every column has at most three vertices.

If a column has three vertices, each of its three row pairs must occur in another column. Splitting the three remaining vertices among two or more columns would supply at most one row pair. They must therefore form a single column with the same three row coordinates. The set $S$ is therefore a $2\times 3$ rectangle. By the $K_3$-uniqueness of $d$, one of its three row pairs has a $\mathcal{D}$-color which occurs on no other row pair. The positive diagonal on this row pair is then the only edge of its amalgamation color, a contradiction.

We may consequently assume that every nonempty column contains one or two vertices. The support of each two-vertex column must be repeated in another such column. There cannot be exactly one such column. If there are exactly three, their row pairs must all coincide, since otherwise one of the pairs would occur only once. In this case $S$ is a $3\times 2$ rectangle. A uniquely colored edge of $c$ on its three columns again gives a uniquely colored positive diagonal. If there are no two-vertex columns, all six column coordinates are distinct, and the $K_6$-uniqueness of $c$ immediately gives a unique amalgamation color.

It remains to consider exactly two two-vertex columns. They have the same row pair, and the other two vertices lie in distinct singleton columns. Write
$$
S=\bigl(\{a,b\}\times\{u,v\}\bigr)
       \cup\{(z,\rho),(w,\sigma)\},
$$
where $a,b,z,w$ are pairwise distinct and $u\ne v$. We impose no distinctness assumption on $\rho,\sigma$ at this stage.

Choose an edge whose $\mathcal{C}$-color is unique on the four columns $\{a,b,z,w\}$. This edge cannot be $ab$, since its two rectangle diagonals have different signs and would each have a unique amalgamation color. It cannot be $zw$, which supports only one edge of $S$. By interchanging the two rectangle columns and the two singleton vertices if necessary, we may therefore suppose it is $az$.

The only two edges of $S$ with this $\mathcal{C}$-color join $(z,\rho)$ to
$(a,u)$ and $(a,v)$. Neither may have a unique amalgamation color, so their full colors must coincide. In particular,
$
\rho\notin\{u,v\}\mbox{ and }  d(\rho u)=d(\rho v)=\delta.
$
Set $\alpha=c(ab)$ and $\beta=d(uv)$.
Since $d$ has no monochromatic triangle, $\delta\ne\beta$.

The rectangle on columns $a,b$ has two diagonals of colors
 $(D,\alpha,\beta,+)$ and $(D,\alpha,\beta,-)$.
Each of these colors must occur on an additional edge of $S$. None of the edges from $(z,\rho)$ to the rectangle can provide such an occurrence, since their $\mathcal{D}$-color is $\delta$.

We claim that at most one edge from $(w,\sigma)$ to the rectangle can have $\mathcal{C}$-color $\alpha$ and $\mathcal{D}$-color $\beta$. If $\sigma\notin\{u,v\}$, then at most one of $d(\sigma u),d(\sigma v)$ equals $\beta$, since otherwise $\{\sigma,u,v\}$ would be a monochromatic triangle in $d$. If $\sigma\in\{u,v\}$, only the other row can supply a diagonal with $\mathcal{D}$-color $\beta$. In either case there are at most two eligible edges, one to each rectangle column. At most one of them has $\mathcal{C}$-color $\alpha$, otherwise $wab$ would be a monochromatic triangle in $c$. This proves the claim.

Consequently one of the two rectangle diagonal colors must be repeated by an edge from $(w,\sigma)$ to the rectangle, and the other must be repeated by the edge $(z,\rho)(w,\sigma)$. Thus
$$
c(zw)=\alpha\mbox{ and }  d(\rho\sigma)=\beta,
$$
and exactly one of $c(wa),c(wb)$ equals $\alpha$. In particular, $c(wa)\ne c(wb)$. The edge joining the two singleton vertices is a diagonal, so $\sigma\ne\rho$. Moreover, since $d(\rho u)=d(\rho v)=\delta\ne\beta$, we have $\sigma\notin\{u,v\}$. The four rows $u,v,\rho,\sigma$ are therefore pairwise distinct.

The repeated diagonal supplied by $(w,\sigma)$ shows that one of $d(\sigma u),d(\sigma v)$ equals $\beta$. Interchanging $u,v$ if necessary, write
$$
d(\sigma u)=\beta\mbox{ and }  d(\sigma v)=\varepsilon.
$$
Here $\varepsilon\ne\beta$, again by the absence of a monochromatic triangle. On the four rows $\{u,v,\rho,\sigma\}$ the colors are now
$$
\begin{array}{c|c}
 \text{row edges}&\text{$\mathcal{D}$-color}\\ \hline
 uv,\ \rho\sigma,\ \sigma u&\beta\\
 \rho u,\ \rho v&\delta\\
 \sigma v&\varepsilon
 \end{array}
$$
If $\varepsilon=\delta$, this $K_4$ has two color classes of size three, contrary to the $K_4$-uniqueness of $d$. Hence $\varepsilon\notin\{\beta,\delta\}$, and $\sigma v$ is the unique row pair with $\mathcal{D}$-color $\varepsilon$.

Exactly two diagonal edges of $S$ have this $\mathcal{D}$-color: the edges from $(w,\sigma)$ to $(a,v)$ and $(b,v)$. Their $\mathcal{C}$-colors are distinct, as $c(wa)\ne c(wb)$. Each therefore has a unique amalgamation color, which is the required contradiction.
\end{proof}

\begin{proof}[\bf Proof of Theorem~\ref{thm:k12}]
Take
$$
\mathcal P_{12}=\{K_3\uniq,K_4\uniq,K_6\uniq,
                        K_8\odd,K_{12}\odd\}.
$$
These conditions hold on one vertex and are preserved under restriction. Let $c$ have property $\mathcal P_{12}$ and let $d$
have property $\U$. Proposition~\ref{prop:unique} preserves
$K_3$- and $K_4$-uniqueness, and
Proposition~\ref{prop:unique6} preserves $K_6$-uniqueness.
Propositions~\ref{prop:eight} and~\ref{prop:twelve} preserve
$K_8$- and $K_{12}$-oddness, respectively. Their hypotheses all follow from $\mathcal P_{12}$ and $\U$.
Thus $c\otimes d$ has property $\mathcal P_{12}$.
Applying Lemma~\ref{lem:iteration} gives the required bound for
$r_{K_{12}}(n)$ and proves Theorem~\ref{thm:k12}.
The construction also retains $K_3$-, $K_4$- and
$K_6$-uniqueness and $K_8$-oddness.
\end{proof}




\end{document}